\documentclass[10.5pt, reqno]{amsart}
\usepackage[margin=3cm]{geometry}

\usepackage{array,booktabs,tabularx}
\usepackage{graphicx}
\usepackage{amssymb,amsmath, amsthm}
\usepackage{enumerate}
\usepackage{xfrac} 
\usepackage{nicefrac}
\usepackage{color}
\usepackage{hyperref}
\usepackage{esint}
\usepackage{mathtools}
\usepackage{dsfont}
\usepackage{ esint }
\usepackage[show]{ed}
\usepackage{comment}
\usepackage{enumitem}
\usepackage{soul}
\usepackage{needspace}
\usepackage[normalem]{ulem}
\usepackage{tikz}
\usepackage{float}

\theoremstyle{plain}

\newtheorem{thm}{Theorem}[section]
\newtheorem{lem}{Lemma}[section]

\newtheorem{prop}{Proposition}[section]
\newtheorem{cor}{Corollary}[section]	
\newtheorem{rem}{Remark}[section]

\allowdisplaybreaks[2]

\newenvironment{proof1}[1]{\begin{trivlist} \item[] {\em Proof of #1:}}{\hfill $\Box$
                      \end{trivlist}}

\newcommand{\ud}{\,\mathrm{d}}

\newcommand{\R}{\mathbb{R}}

\newcommand{\Om}{\Omega}

\newcommand{\pa}{\partial}

\newcommand{\norm}[1]{\left\lVert#1\right\rVert}
\newcommand{\tih}{\tilde{h}}
\newcommand{\tE}{\tilde{E}}

\newcommand{\tS}{\tilde{S}} 

\newcommand{\pb}{\phi_{_{\!B}}\!}  
\newcommand{\pt}{\phi_{_{\!T}}\!} 
\newcommand{\pl}{\phi_{_{\!L}}\!}   
\newcommand{\pr}{\phi_{_{\!R}}\!}  
\newcommand{\rhb}{\rho_{_{\!B}}\!}  
\newcommand{\rht}{\rho_{_{\!T}}\!} 
\newcommand{\rhl}{\rho_{_{\!L}}\!}   
\newcommand{\rhr}{\rho_{_{\!R}}\!}  

\newcommand{\x}{x_{{0}}}  
\newcommand{\y}{y_{{0}}}  
\newcommand{\xx}{x_{{1}}}  
\newcommand{\yy}{y_{{1}}}

\DeclareMathOperator{\diam}{diam}

\title{Nodal line estimates for the second Dirichlet eigenfunction}

\author[T. Beck]{Thomas Beck}
\email{tdbeck@email.unc.edu}
\address{Department of Mathematics, University of North Carolina at Chapel Hill \\ CB\#3250
  Phillips Hall \\ Chapel Hill, NC 27599}
  
  \author[Y. Canzani]{Yaiza Canzani}
\email{canzani@email.unc.edu}
\address{Department of Mathematics, University of North Carolina at Chapel Hill \\ CB\#3250
  Phillips Hall \\ Chapel Hill, NC 27599}

\author[J.L. Marzuola]{Jeremy L. Marzuola}
\email{marzuola@math.unc.edu}
\address{Department of Mathematics, University of North Carolina at Chapel Hill \\ CB\#3250
  Phillips Hall \\ Chapel Hill, NC 27599}

\date{\today}                                        

\begin{document}
\maketitle

\begin{abstract}
 We study the nodal curves of low energy Dirichlet eigenfunctions in generalized curvilinear quadrilaterals.  The techniques can be seen as a generalization of the tools developed by Grieser-Jerison in a series of works on convex planar domains and rectangles with one curved edge and a large aspect ratio.  Here, we study the structure of the nodal curve in greater detail, in that we find precise bounds on its curvature, with uniform estimates up to the two points where it meets the domain at right angles, and show that many of our results hold for relatively small aspect ratios of the side lengths.  We also discuss applications of our results to Courant-sharp eigenfunctions and spectral partitioning.
\end{abstract}

%%%%%%%%%%%%%%%%%%%%%%%%%%%%%%%%%%%%%
%%%%%%%%%%%%%%%%%%%%%%%%%%%%%%%%%%%%%
\section{Introduction and statement of results}
%%%%%%%%%%%%%%%%%%%%%%%%%%%%%%%%%%%%%
%%%%%%%%%%%%%%%%%%%%%%%%%%%%%%%%%%%%%
%\begin{itemize}
%    \item Importance of low lying eigenstates;
%    \item Meaning of the zero set for low lying eigenfunctions;
%    \item G-J and our improvements;
%    \item Naive partitioning;
%    \item Results
%\end{itemize}

Understanding the fundamental modes of vibration of a compact domain is a longstanding problem. The original motivation was to describe how a metal sheet with a given shape would vibrate when struck at some fundamental frequency. The main goal being to understand the structure of the set of points in the sheet that are stable, i.e. that are not vibrating.  These non-vibrating regions are the zero sets of the Laplace eigenfunctions corresponding to solving the Helmholtz equation on the domain that represents the metal sheet. In the 17th century R. Hook observed these patterns by spilling sand on a glass sheet, and striking the sheet with a violin bow. When the sheet starts vibrating, the sand rearranges itself across the sheet until it is placed on the non-vibrating areas, thus exhibiting the zero sets for the corresponding eigenfunction. This experiment was later reproduced by E. Chladni, who was the first to record an extensive list of zero set configurations. It is nowadays known as the Chladni plates experiment. 

We dedicate this article to giving a precise description of the structure of the zero set of the second eigenfunction for a planar domain whose shape is obtained after perturbing a rectangle.  While we focus on the second Dirichlet eigenfunction, the techniques developed here can be applied to the low-lying eigenfunctions in general up to a frequency depending upon the length of the domain.  Also, there are natural generalizations to Neumann (or more generally Robin) boundary conditions, but for the sake of clarity and presentation we will focus on Dirichlet domains at present.

We note that low-lying eigenvalues and eigenfunctions of the Laplacian on a compact domain also play a role in understanding random walks (\cite{kenig1989h}), heat conductivity (\cite{smits1996spectral}) and more. See for instance the recent works \cite{zelditch2017eigenfunctions,steinerberger2017localization} and references therein for a nice overview of applications and modern topics in the theory of eigenfunctions and nodal sets.

%Low-lying eigenvalues and eigenfunctions of the Laplacian on a compact domain play a role in understanding random walks (\cite{kenig1989h}), fundamental modes of vibration (\cite{kac1966can}), heat conductivity (\cite{smits1996spectral}) and more.  The nodal sets of eigenfunctions represent regions of the domain that do not vibrate at a given frequency.  In particular, the low lying eigenmodes are the lowest energy modes for a quantum state confined to a given geometry satisfying the linear Schr\"odinger equation $i \partial_t \psi (x,t) + \Delta \psi(x,t) = 0$, $\psi(x,0) \in L^2 (\Omega)$, $\psi(t,\cdot)|_{\partial \Omega} = 0$. 
%See for instance the recent works \cite{zelditch2017eigenfunctions,steinerberger2017localization} and references therein for a nice overview of applications and modern topics in the theory of eigenfunctions and nodal sets.

In this work, we study the second Dirichlet eigenfunction of the Laplacian on a planar domain $\Omega$, so that 
\begin{eqnarray*}
    \left\{ \begin{array}{rlc}
    \Delta v(x,y) & = -\mu v(x,y) & \text{in } \Om\\
   v(x,y) & = 0& \text{on } \pa \Om
    \end{array} \right.
\end{eqnarray*}
where $\mu$ is the corresponding eigenvalue. The domain $\Omega$ is a curvilinear rectangle that is very nearly rectangular in a Gromov-Hausdorff sense to be made precise below in \eqref{eqn:Omega}. For convenience, we normalize $v$ so that $\norm{v}_{L^{\infty}} = 1$. We are interested in studying the nodal set of $v$, which we denote by
\begin{align*}
\Gamma = \{(x,y)\in\mathring{\Omega}:v(x,y) = 0\}. 
\end{align*}
 To do this we build off the pioneering works of Jerison \cite{jerison-nodal} and Grieser-Jerison \cite{grieser1996asymptotics,grieser-maximum,GJ-Pacific}, who studied the low energy eigenfunctions in convex domains and rectangles of high aspect ratio with one curved edge. For convex domains they studied the location of the maximum and nodal line of the first and second Dirichlet eigenfunction respectively, giving estimates that are uniform as the eccentricity of the domain increases. They also derived a method to do a very detailed asymptotic analysis of the location of the nodal line or the location of the maximum for low energy Dirichlet eigenfunction in a rectangle with one curved side. This method is the starting point for our work on curvilinear rectangles.
 
On a rectangle, the zero set of the second eigenfunction for the Laplacian is a straight line, perpendicular to the long sides, that divides the rectangle in two equal pieces. Here, we study the zero set of the second eigenfunction $v$ on a region $\Omega$ that is a perturbation of a rectangle. Indeed, we extend the results in \cite{GJ-Pacific} to explore the precise dependence of the nodal set 
 on the properties of the bounding curves and the aspect ratio of the underlying  region $\Omega$. 
 
 We obtain estimates on the width and regularity of the nodal line of $v$, with explicit bounds, tracking how the slope and curvature of the nodal line depend upon the top and bottom curve. In particular, we show that there are distinct differences in the nodal line depending upon if some of the bounding curves are flat or if they are curved. As we are imposing Dirichlet boundary conditions, the eigenfunction $v$ vanishes on the boundary of the domain, and analyzing the behavior of the nodal line becomes increasingly delicate as it approaches the boundary. Our techniques allow us to obtain estimates that are uniform up to the boundary and show that the nodal line meets the boundary of $\Omega$ orthogonally at two points. 
 
%  In  addition, with the intention of motivating the iteration of our analysis to explore repeated nodal cuts, we study in detail the behavior of the nodal line as it approaches the boundary as well as bounds on the curvature of the nodal line within the region.  In addition, we show distinct and strong differences in the nodal line depending upon if some of the bounding curves are flat or if they are curved.  In the case they are curved, our methods give tools for tracking the dependence upon the top and bottom curve.

To state our results precisely, we first define the class of domains $\Omega$ under consideration. Let $\pb, \pt, \pl, \pr$ be functions defining a region $\Omega$ in $\R^2$ that is a perturbation of the rectangle $[0, N]\times [0,1]$,  for $N>0$, of the form
\begin{align} \label{eqn:Omega}
\Omega= \{(x,y)\in\R^2:\;\;  \pb(x) \leq y \leq \pt(x),\;\;  \pl(y) \leq x \leq \pr(y)\}.
\end{align}
\begin{center}
    \includegraphics[width=9cm]{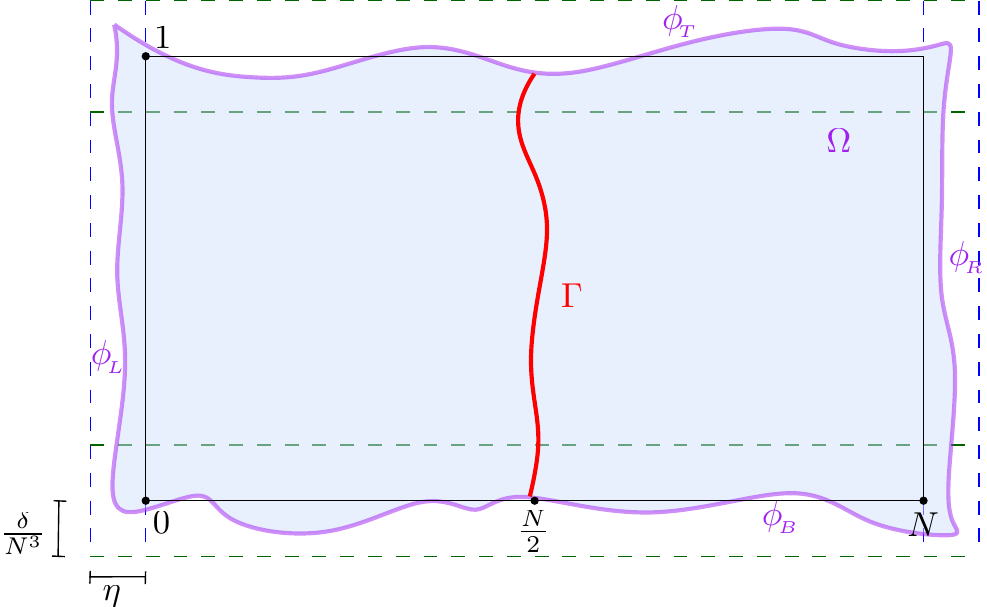}
\end{center}
The functions $\pl$, $\pr$ defining the sides of the domain are in $C^{2}([-\tfrac{1}{2},\tfrac{3}{2}];\mathbb{R})$, with
\begin{align} \label{eqn:pl}
-\eta\leq \pl\leq 0,\qquad 
N \leq \pr \leq \eta+N, \qquad 
\left|\frac{d^j}{dy^j}\pl\right|\leq \eta, \qquad
\left|\frac{d^j}{dy^j}\pr\right| \leq \eta,
\end{align}
for $j=1,2$ and some $\eta>0$. The functions $\pb$, $\pt$ defining the top and bottom are in $C^{\infty}([-\tfrac{1}{2},N+\tfrac{1}{2}];\mathbb{R})$, with
\begin{align} \label{eqn:pb}
\left|\pb\right| \leq \frac{\delta}{N^3}, \qquad
\left|\pt - 1\right|\leq \frac{\delta}{N^3}, \qquad 
\left|\frac{d^j}{dx^j}\pb\right|\leq \tilde{C}_j\frac{\delta}{N^3},\qquad  \left|\frac{d^j}{dx^j}\pt\right| \leq \tilde{C}_j\frac{\delta}{N^3},
\end{align}
for $j\geq 1$ and constants $\tilde{C}_j>0$. Here  $\eta>0,$ $0<\delta<\tfrac{1}{5}$, and $N\geq 5$. Note that as $\eta$, $\delta$ tend to $0$, the domain $\Omega$ becomes rectangular. 

In the case of the rectangle $[0,N]\times[0,1]$ with $N>1$, the second Dirichlet eigenfunction is given by $v(x,y) = \sin\left(\tfrac{2\pi x}{N}\right)\sin\left(\pi y\right)$ and the nodal line $\Gamma$ is precisely the straight line $\tfrac{N}{2}\times(0,1)$. The theorem below shows how $\Gamma$ changes under the above perturbations of the rectangle.
Let $\pi_x:\R^2 \to \R$ be the projection onto the $x$-axis.

%%%%%%%%%%%%%%%%%%%%%%%%%%%%%%%%%%%%%

\begin{thm} \label{thm:Main}
There exist  $c>0$, $C>0$, such that  $\pi_x(\Gamma)\subset [\tfrac{N}{2}-CN(\eta  + \tfrac{\delta}{N}),\tfrac{N}{2}+CN(\eta  + \tfrac{\delta}{N})]$ and has diameter bounded by
\begin{align*}
\diam(\pi_x(\Gamma)) \leq C\left(\eta e^{-c N} + \frac{\delta}{N^2}\right).
\end{align*}
Moreover, there exists a function $g(y)$ such that $\Gamma\cap  \mathring{\Omega}= \{(x,y)\in\mathring{\Omega}: x = g(y)\}$, with 
\begin{align*}
\left|g'(y)\right| +\left|g''(y)\right| \leq C\left(\eta e^{-c N} + \frac{\delta}{N^2}\right).
\end{align*}
The nodal line $\Gamma\cap  \mathring{\Omega}$ touches the boundary of $\Omega$ at precisely $2$ points, and it meets the boundary orthogonally at these points.
\end{thm}
%%%%%%%%%%%%%%%%%%%%%%%%%%%%%%%%%%%%%
Here and throughout, constants denoted by $c$, $C$, $C_1$ etc depend on the constants $\tilde{C}_j$, but are independent of $\eta$, $\delta$, and $N$. (In fact we will only require control on derivatives up to $j=5$.) 

An immediate feature to note is that in the special case of flat upper and lower boundaries ($\pb(x) \equiv 0$, $\pt(x)\equiv1$), we can set $\delta = 0$ and the factor of $\tfrac{\delta}{N^2}$ does not appear in the estimates of Theorem \ref{thm:Main}. Therefore, in this flat case the diameter of the nodal line, $\diam(\pi_x(\Gamma))$, is exponentially small in $N$ (rather than the polynomial decay in $N$ when $\delta\neq0$).

From Theorem \ref{thm:Main} we see that for $N$ sufficiently large (and $\delta \leq C \eta$), the perturbation of the nodal line from straight is smaller than that of the side perturbations $\pl(x)$, $\pr(x)$. In the flat case, $\pb(x) \equiv 0$, $\pt(x)\equiv1$, we can track the constants in the proof of Theorem \ref{thm:Main} (see Section \ref{Fodes}) to obtain an explicit lower bound on the size of $N$ required for this to occur:
%%%%%%%%%%%%%%%%%%%%%%%%%%%%%%%%%%%%%
\begin{cor} \label{cor:Main}
There exists a constant $N_0>0$ such that for $N\geq N_0$ and $\delta\leq\eta$, 
\begin{align*}
\diam(\pi_x(\Gamma)) \leq  \tfrac{\eta}{2}\qquad\text{and}\qquad   \left|g'(y)\right| +\left|g''(y)\right| \leq \tfrac{\eta}{2}.
\end{align*} 
In the flat case, $\pb(x) \equiv 0$, $\pt(x)\equiv1$, for each $N\geq8$, we can take $\eta = \eta(N)>0$ sufficiently small so that the above estimates hold. 
\end{cor}
%%%%%%%%%%%%%%%%%%%%%%%%%%%%%%%%%%%%%

%\begin{rem} \label{rem:flat}
%In the special case where $\pb(x) \equiv 0$, $\pt(x)\equiv1$, the factor of $\tfrac{\delta}{N^2}$ is not required in the statement of Theorem \ref{thm:Main}, and the width of the nodal line is exponentially small in $N$ (rather than polynomial). In Section \ref{Fodes} we will show that by tracking the constants in the proof of Theorem \ref{thm:Main}, for each $N\geq8$, we can take $\eta = \eta(N)>0$ sufficiently small so that the estimates in Corollary \ref{cor:Main} hold. %({\color{red} We haven't yet tracked the size of $N$ needed to get the $g''(y)$ bound in the flat case})
%\end{rem}
By controlling the behavior of the nodal line up to the boundary, we are able to show for the class of domains under consideration that the nodal line is not closed, but meets the boundary (orthogonally) at two points. 
%\begin{rem} \label{rem:orthogonal}
 More generally, Payne \cite{Payne} conjectured that the nodal line of the second eigenfunction of a bounded planar domain touches the boundary at $2$ points. This was proved for smooth, convex domains by Melas \cite{Melas}, but a counterexample (for a non-simply connected planar domain) was given by Hoffmann-Ostenhof, Hoffmann-Ostenhof, and Nadirashvili \cite{HON}.
 
In \cite{freitas-krej}, Freitas and Krej\u{c}i\u{r}\'ik study the Dirichlet Laplacian for a class of thin curved tubes. As the volume of the cross-section tends to $0$ they establish the convergence of the eigenvalues and eigenfunctions in terms of an ordinary differential operator on the base curve of the tube. In particular, they locate the nodal set to sufficient precision to also deduce that the nodal set must intersect the boundary. Krej\u{c}i\u{r}\'ik and Tu\u{s}ek also prove an analogous result for domains consisting of a thin tubular neighborhood of a hypersurface, \cite{krej-tusek}. The idea of reducing to an associated ordinary differential operator has also been used extensively by Friedlander-Solomyak \cite{fried-sol1}, \cite{fried-sol2} and Borisov-Freitas \cite{bor-frei} to obtain asymptotics of the eigenvalues, eigenfunctions, and resolvent of the Dirichlet Laplacian in thin domains.

\subsection*{Applications to partitioning algorithms}
Recently, in work on graph and data partitioning algorithms, Szlam et al in \cite{szlam1,szlam2} observed that if one partitions very general geometric graphs using cuts along nodal lines of the graph Laplacian, that the regions tend towards rectangles of bounded aspect ratios.  The underlying idea of graph partitions are for instance to cluster data points or to provide a good foundation for a wavelet basis to name just two.  There is also the continuum limit version of this, in which one could ask to partition a planar domain using the first non-trivial Neumann or second Dirichlet eigenfunction respectively.  Again, for a very general boundary, one expects that such partitions would converge rapidly to a set of near rectangles.  

We can use Corollary \ref{cor:Main} to demonstrate such a convergence: In the flat case, with $N\geq8$, and $\eta =\eta(N)$ sufficiently small, we perform the following iteration procedure. Given such a domain $\Omega$, with flat top and bottom boundaries, we form $2$ new domains by \textit{cutting} $\Omega$ along the nodal line $\Gamma$. By then rescaling each domain in the $x$-direction, and using the estimates above, we can ensure that these two new domains are of the form of $\Omega$, for appropriate chosen functions $\pl^{\!new}(y)$ and $\pr^{\!new}(y)$, and satisfy the bounds
\begin{align*}
    \tfrac{3}{4}\eta\leq \pl^{\!new}\leq 0, \quad 0 \leq \pr^{\!new} - N \leq \tfrac{3}{4}\eta, \quad \left|\frac{d^j}{dy^j}\pl^{\!new}\right|, \left|\frac{d^j}{dy^j}\pr^{\!new}\right| \leq \frac{3}{4}\eta.
\end{align*}
In other words, these new sides to the domain satisfy the original bounds in \eqref{eqn:pl}, but with $\eta$ replaced by $\tfrac{3}{4}\eta$. Iterating this process of cutting along the nodal line, will therefore give domains converging (in a Hausdorff sense) to an exact rectangle. Analogously, for the general top and bottom boundaries considered here, using the estimates in Corollary \ref{cor:Main}, given $\eta$, $\delta$, with $\delta \leq \eta$, and for $N$ sufficiently large, we can repeat the above iteration process, to again give a sequence of domains converging to a rectangle.

Another partitioning related to the nodal set of Dirichlet eigenfunctions of the Laplacian is the following: Given a domain $\Omega$ and integer $k\geq2$, a spectral minimal k-partition of $\Omega$ is a partition of $\Omega$ into $k$ disjoint sets $\Omega_i$ that minimizes $\max_{i}\lambda(\Omega_i)$. Here $\lambda(\Omega_i)$ is the first Dirichlet eigenvalue of $\Omega_i$. If $k=2$, then the spectral minimal partition is given by the nodal domains of a second Dirichlet eigenfunction of $\Omega$. More generally, if a $k$-th Dirichlet eigenfunction has exactly $k$-nodal domains (and so gives equality in the Courant nodal domain Theorem), then these nodal domains form a minimal $k$-partition. See the survey paper of Helffer \cite{helffer-minimal1} for greater discussion of spectral minimal partitions and references. It is therefore important to classify examples where the Courant nodal domain Theorem is sharp, and in \cite{helffer-minimal2} they show that the third Dirichlet eigenfunction of the rectangle has three nodal domains whenever the aspect ratio is greater than $\sqrt{8/3}$. Using the techniques presented here, for any fixed $j\geq2$, by taking $N$ sufficiently large, for all $2\leq k \leq j$, the $k$-th Dirichlet eigenfunction has exactly $k$ nodal domains (with nodal set approximately equal to the union of the $k-1$ lines $\{\tfrac{l}{kN}\}\times[0,1]$ for $1\leq l \leq k-1$). Thus, in this case, the nodal domains will provide a spectral minimal $k$-partition.

\subsection*{Outline of the paper}
The structure of the rest of the paper as follows: In Section \ref{sec:Adiabatic} we describe an adiabatic approximation of the eigenfunction that is a key ingredient in the proof of Theorem \ref{thm:Main}. This type of approximation, which can be viewed as an approximate separation of variables for our approximately rectangular domain, has been used in the work of Grieser and Jerison \cite{grieser1996asymptotics}, \cite{GJ-Pacific}. The approximation has also been used in \cite{BSS} for numerical analysis of eigenfunctions in partially rectangular billiards, and in \cite{HM1} to analyze non-concentration of eigenfunctions in partially rectangular billiards. In Section \ref{sec:nodal}, we establish the desired properties of the width and regularity of the nodal line using the adiabatic approximation. Then, in Section \ref{duhamel} we demonstrate how in the flat case, we have simple ODE estimates to establish the approximation, and following this, we prove the error estimates for the approximation for our general class of domains.  Lastly, in Section \ref{Fodes} we compute an explicit Hadamard variation formula to evaluate the effect the side perturbations have on the eigenfunction. This will in particular allow us to track the constants appearing in the proof of Theorem \ref{thm:Main} in the flat case and prove Corollary \ref{cor:Main}.

\subsection*{Acknowledgements} This project was started due to a conversation with Stefan Steinerberger and Hau-tieng Wu, who pointed us to the work of \cite{szlam2} as motivation to understand nodal line partitioning in domains, and with whom the third author is exploring a related question for the $1$-Laplacian on curvilinear rectangles.  YC is supported in part by the Sloan Foundation.  JLM acknowledges support from the NSF through NSF CAREER Grant DMS-1352353.

%%%%%%%%%%%%%%%%%%%%%%%%%%%%%%%%%%%%%
%%%%%%%%%%%%%%%%%%%%%%%%%%%%%%%%%%%%%
\section{The Adiabatic Ansatz} \label{sec:Adiabatic}
%%%%%%%%%%%%%%%%%%%%%%%%%%%%%%%%%%%%%
%%%%%%%%%%%%%%%%%%%%%%%%%%%%%%%%%%%%%
A key ingredient in the proof of Theorem \ref{thm:Main} is to establish properties of a Fourier decomposition of the eigenfunction $v(x,y)$. For convenience, we introduce the height function
\[
h(x):=\pt(x)-\pb(x),
\]
and note that \eqref{eqn:pb} we have $1-\tfrac{2\delta}{N^3} \leq h(x) \leq 1+\tfrac{2\delta}{N^3}$ for all $x\in[0,N]$. For  $(x,y)\in\Omega$ with $0\leq x\leq N$, we write $v$ as
\begin{align} \label{eqn:v-decom}
v(x,y) = v_1(x) \sin(\beta(x,y))+ E(x,y),
\end{align}
where
\[
\beta(x,y):= \frac{\pi (y-\pb(x))}{h(x)}.
\]
We will view the first term in the right hand side of \eqref{eqn:v-decom} as the main term, with $E(x,y)$ an error term when $N$ is large, and $\delta$, $\eta$ are small. The function $v_1(x)$ is the first Fourier mode in the $y$-direction, given by
\begin{align*}
v_1(x) = \frac{2}{h(x)}\int_{\pb(x)}^{\pt(x)} v(x,y)\sin(\beta(x,y)) \ud y.
\end{align*}
To prove Theorem \ref{thm:Main} we will use this decomposition of $v$, and will require a lower bound on $|v_1'(x)|$, together with upper bounds on $v_1(x)$, $E(x,y)$ and their derivatives. In fact, to prove the estimates on regularity of the nodal line near $\pa\Omega$, we need to consider a larger class of decompositions of $v$: Given $(\x,\y)\in \Gamma$, suppose that  $(\xx, \yy)\in \partial \Omega $, with $\yy = \pb(\xx)$ and $d( (\x,\y), \partial \Omega)=  d((\x,\y),(\xx,\yy))$. Next,  let $F$ the linear isometry built by rotating about the point $(\x,\y)$ and then translating in such a way that if $(\x,\y)=F(\tilde{x}_0, \tilde{y}_0)$, then $(\xx,\yy)=F(\tilde{x}_0,0)$. In general, write $(\tilde x, \tilde y)=F^{-1}(x,y)$ for the new system of coordinates. We then define $w(\tilde{x},\tilde{y})$ to be equal to the eigenfunction $v$ in these rotated coordinates,
\[
w=v\circ F.
\]
\begin{rem} \label{rem:rotation}
By the bounds on $\pb'(x)$ and $\pt'(x)$ from \eqref{eqn:pb}, there exists $C>0$ such that the angle of rotation is bounded by $C\tfrac{\delta}{N^3}$. 
\end{rem}
 The function $w$ satisfies \[(\Delta + \mu)w = 0\] in the domain
\begin{align*}
\tilde{\Omega} = \{(\tilde{x},\tilde{y}):\rhl(\tilde{y})\leq \tilde{x} \leq \rhr(\tilde{y}),\; \rhb(\tilde{x}) \leq \tilde{y} \leq \rht(\tilde{x})\},
\end{align*}
with $w|_{\pa\tilde{\Omega}} = 0$. In particular, for $0\leq \tilde{x} \leq N$,  we have $w(\tilde{x},\rhb(\tilde{x})) = w(\tilde{x},\rht(\tilde{x})) = 0$. Here $\rhb$, $\rht$ satisfy the bounds
\begin{align*}
\left|\rhb(\tilde{x})\right| \leq \frac{2\delta}{N^3}(1+|\tilde{x}-\tilde{x}_0|), \qquad\qquad \left|\rht(\tilde{x}) - 1\right|\leq \frac{2\delta}{N^3}(1+|\tilde{x}-\tilde{x}_0|), 
\end{align*}
\begin{align}\label{e:rho B T}
\left|\frac{d^j}{d\tilde{x}^j}\rhb(\tilde{x})\right|\leq 2\tilde{C}_j\frac{\delta}{N^3},\qquad \qquad\qquad \left|\frac{d^j}{d\tilde{x}^j}\rht(\tilde{x})\right| \leq 2\tilde{C}_j\frac{\delta}{N^3},
\end{align}
$j\geq1$. Up to the factor of $2$, the derivative bounds are the same as for $\pb$, $\pt$. Moreover, by the construction of $F$, we have $\rhb(\tilde{x}_0) = \rhb'(\tilde{x}_0)=0$. The functions $\rho_L(\tilde{y})$, $\rho_{R}(\tilde{y})$ satisfy
\begin{align*}
    -\eta- \frac{\delta}{N^3} \leq \rhl(\tilde{y})\leq \frac{\delta}{N^3}, \qquad\qquad -\frac{\delta}{N^3} \leq \rhr(\tilde{y}) - N \leq \eta + \frac{\delta}{N^3},
\end{align*}
\begin{align}\label{e:rho L R}
\left|\frac{d^j}{d\tilde{y}^j}\rhl(\tilde{y})\right| \leq \eta + \frac{\delta}{N^3}, \qquad\qquad \left|\frac{d^j}{d\tilde{y}^j}\rhr(\tilde{y})\right| \leq \eta + \frac{\delta}{N^3},
\end{align}
for $j=1,2$. We can make the analogous definition if the closest point to $(\x,\y)$ lies on the upper boundary of $\Omega$. For ease of notation, we now drop the tildes, and for each function $w(x,y)$ coming from such a rotation, for $x\in[0,N]$ we write
\begin{align} \label{eqn:w-decom}
w(x,y) = w_1(x) \sin(\tilde{\beta}(x,y))+ \tilde{E}(x,y),
\end{align}
where
\[
\tilde{\beta}(x,y):= \frac{\pi (y-\rhb(x))}{\tilde{h}(x)},
\]
for the new height function $\tilde{h}(x) = \rht(x)-\rhb(x)$.  To prove Theorem \ref{thm:Main}, we will use the proposition below which gives properties of these decompositions.
\begin{prop} \label{prop:Main}
There exist positive constants $c$, $C$ such that the following properties hold: For each decomposition, there exists a unique point $x_0\in[\tfrac{1}{4}N,\tfrac{3}{4}N]$ such that $w_1(x_0)=0$, and this point lies in the interval $[\tfrac{N}{2}-CN(\eta +\delta),\tfrac{N}{2}+CN(\eta +\delta)]$. Moreover, for $x\in[\tfrac{1}{4}N,\tfrac{3}{4}N]$, \[|w_1'(x)|\geq C^{-1}N^{-1},\]
and for $x\in[1,N-1]$, $0\leq j\leq3$, we have
\begin{align*}
  |w_1^{(j)}(x)| \leq CN^{-j}, \qquad \sup_{y\in[\rhb(x),\rht(x)]}\left|\nabla^j\tilde{E}(x,y)\right| \leq C\left(\eta e^{-cN} + \frac{\delta}{N^3}\right). 
\end{align*} 
\end{prop}
Proposition \ref{prop:Main} is proved in Section \ref{sec: main prop}.
\begin{rem} \label{rem:trivial}
When the rotation is trivial, $w$ is equal to $v$ and the decomposition reduces to the one for $v$ given in \eqref{eqn:v-decom}. Therefore, the properties in this proposition also hold for $v_1(x)$ and $E(x,y)$. In fact, in this case, the unique point $x_0$ where $v_1(x_0)=0$ lies in the interval \[[\tfrac{N}{2}-CN(\eta +\tfrac{\delta}{N})\;,\;\tfrac{N}{2}+CN(\eta +\tfrac{\delta}{N})].\]
\end{rem}

%%%%%%%%%%%%%%%%%%%%%%%%%%%%%%%%%%%%%
%%%%%%%%%%%%%%%%%%%%%%%%%%%%%%%%%%%%%
\section{Estimates on the nodal line} \label{sec:nodal}
%%%%%%%%%%%%%%%%%%%%%%%%%%%%%%%%%%%%%
%%%%%%%%%%%%%%%%%%%%%%%%%%%%%%%%%%%%%
In this section we will prove Theorem \ref{thm:Main} assuming that Proposition \ref{prop:Main} holds. We first establish an upper bound on the width of the projection to the $x$-axis of the nodal line in terms of the error $E$ and its derivatives. We will require a different argument to  control the behavior of the nodal line near the boundary, and so we set \[S(x)=S^{\!B}(x)\cup S^{\!T}(x):=[\pb(x), \pb(x)+\tfrac{1}{4}] \cup [\pt(x)-\tfrac{1}{4}\,,\, \pt(x)].\] 

We continue to write $h(x)=\pt(x)-\pb(x)$ and $\beta(x,y)=\pi(y-\pb(x))/h(x)$.
Since $h(x)\geq \tfrac{1}{2}$ for all $x$, we have that $0 \leq \beta(x,y) \leq \frac{\pi}{2} $ on $S^{\!B}$ and
$\tfrac{3 \pi}{2} \leq \beta(x,y) \leq \pi $ on $S^{\!T}$. Therefore, this choice 
yields  
\[\sin(\beta(x,y)) \geq \frac{2}{\pi} \beta(x,y) =  \frac{2 (y-\pb(x))}{h(x)} \quad  y \in S^{\!B}(x), \]
and similarly
\[
\sin(\beta(x,y)) \geq \frac{ 2 (\pt(x)-y)}{h(x)} \quad  y \in S^{\!T}(x).\] 

Using the decomposition of $v$ from \eqref{eqn:v-decom}, define $\tilde I$ to be the smallest interval with $\x \in \tilde I$  and such that
\begin{enumerate}
\item[A)] $\displaystyle{\sup_{\substack{x \in [0,N]\\ y\in S(x)^c}}|E(x,y)|<  \frac{1}{2}  \inf_{x\in \tilde I^c}\frac{|v_1(x)|}{h(x)}}$
\item[B)] $\displaystyle{\sup_{\substack{x \in [0,N]\\ y \in S(x)}}|\partial_y E(x,y)|< 2 \,  \inf_{x\in \tilde I^c}\frac{|v_1(x)|}{h(x)}}$
\end{enumerate}

%%%%%%%%%%%%%%%%%%%%%%%%%%%%%%%%%%%%%%%%%%%%%%%%%%%%%

\begin{lem}\label{lem:interval}
If $x\in  \tilde I^c$, then $v(x,y)\neq 0$ for all $y \in (\pb(x),\pt(x))$.
\end{lem}
\begin{proof1}{Lemma \ref{lem:interval}}
Let $y \in S(x)^c$. Then,   $ \sin(\beta(x,y))\geq \sin(\beta(x,\pb(x)+\tfrac{1}{4}))\geq  \frac{1}{2 h(x)}$ and so 
\begin{align} \label{eqn:interval1}
|v(x,y)|= | v_1(x)  \sin(\beta(x,y))+ E(x,y)| \geq  |v_1(x)|   \tfrac{1}{2h(x)}  - |E(x,y)|>0.
\end{align}

By assumption (A), this is strictly positive. Now let $y \in S^{\!B}(x)$. Then,  since $E(x, \pb(x))=0$, we have $|E(x,y)| \leq (y-\pb(x)) \sup_{y \in S^{\!B}(x)}|\partial_yE(x,y)| $.  Also, using that $\sin(\beta(x,y))\geq { \tfrac{2(y-\pb(x))}{h(x)}}$  we obtain
\begin{align} \label{eqn:interval2}
|v(x,y)|= | v_1(x)\sin(\beta(x,y)) + E(x,y)| \geq (y-\pb(x))\Big( |v_1(x)|  { \tfrac{2 }{h(x)}} -  \sup_{y \in S^{\!B}(x)}|\partial_y E(x,y)| \Big)>0,
\end{align}
and by assumption (B) this is strictly positive. The case $y \in S^{\!T}(x)$ is treated in the same way.
\end{proof1}

Using Proposition \ref{prop:Main}, we let $x_0$ be the unique point in the interval $[\tfrac{1}{4}N,\tfrac{3}{4}N]$ where $v_1(x_0)=0$. Define the interval $I$,
\[
I := [x_0-\tau,x_0+\tau]%I:=[\x- \tau_{_{\!\phi\!, N\!, E,\alpha}}\;,\; \x +\tau_{_{\!\phi\!, N\!, E,\alpha}}],
\]
where
\[
\tau :=\frac{\sup_{x\in\tilde{I}}h(x)}{ 2 \inf_{x\in \tilde I}|v_1'(x)|} \max\Big\{ 4\sup_{\substack{(x,y)\in {\Omega} \\ x \in \tilde{I}}}|E(x,y)|\,,\,   \sup_{\substack{(x,y)\in {\Omega} \\ x \in \tilde{I}}}|\partial_yE(x,y)| \Big\}.
\]

%%%%%%%%%%%%%%%%%%%%%%%%%%%%%%%%%%%%%%%%%%%%%%%%%%%%%
\begin{lem}\label{lem:interval2}
If $x\in  I^c$, then $v(x,y)\neq 0$ for all $y \in (\pb(x),\pt(x))$.
\end{lem}

\begin{proof1}{Lemma \ref{lem:interval2}}
Let $y \in S(x)^c$. Then, as in \eqref{eqn:interval1},
\[
|v(x,y)| \geq  |v_1(x)|  \tfrac{1}{2h(x)}  - |E(x,y)|.
\]
Also, since $v_1(\x)=0$, we have $|v_1(x)|\geq |x-\x|\inf_{x\in \tilde I}|v_1'(x)|$.  Therefore,  $|v(x,y)|>0$ provided
\[
|x-\x| > \frac{ 2 h(x)\, \sup_{y}|E(x,y)|}{ \inf_{x\in \tilde I}|v_1'(x)|   },
\]
and the latter always holds if $|x-\x|> \tau$.

Now let $y \in  S^{\!B}(x)$. Then, as in \eqref{eqn:interval2}
\begin{align*}
|v(x,y)|
\geq (y-\pb(x))\Big( |v_1(x)| \tfrac{2 }{h(x)} -  \sup_{y \in S^{\!B}(x)}|\partial_yE(x,y)| \Big)
\end{align*}
and using $|v_1(x)|\geq |x-\x|\inf_{x\in \tilde I}|v_1'(x)|$ gives
\begin{align*}
|v(x,y)| \geq (y-\pb(x))\Big( \tfrac{2 }{h(x)} |x-\x|\inf_{x\in \tilde I}|v_1'(x)| -  \sup_{y \in S^{\!B}(x)}|\partial_yE(x,y)| \Big)
\end{align*}
Therefore, from the definition of $\tau$,  
$|v(x,y)|>0$ for $|x-x_0|>\tau$.
\end{proof1}

Using Proposition \ref{prop:Main} and Remark \ref{rem:trivial}, there exist   $c>0$ and $C>0$ such that $\tau \leq C\left(\eta e^{-cN}+\tfrac{\delta}{N^2}\right)$
 and $I\subset[\tfrac{N}{2} - CN(\eta + \delta N^{-1}), \tfrac{N}{2} + CN(\eta + \delta N^{-1})]$, and so the estimate on the width and location of the nodal line in Theorem \ref{thm:Main} follows immediately from Lemma \ref{lem:interval2}.

To study the regularity of the nodal line, we use the coordinate change described in Section \ref{sec:Adiabatic}. For a given $(x_0,y_0)$ with $v(x_0,y_0) = 0$, $y_0\leq\tfrac{1}{2}$, this coordinate change transforms $(x_0,y_0)$ to $(\tilde{x}_0,\tilde{y}_0)$ and the eigenfunction $v$ to $w(\tilde{x},\tilde{y})$. Dropping the tildes, we have $(\Delta+\mu)w = 0$ in the domain 
\begin{align}
\label{tOmegadef}
\tilde{\Omega} = \{(x,y):\rho_L(y)\leq x \leq \rho_R(y), \rhb(x) \leq y \leq \rht(x)\},
\end{align}
with $w(x,\rhb(x)) = w(x,\rht(x)) = 0$. Moreover, $\rhb(x_0) = \rhb'(x_0) = 0$. Setting $\tilde{h}(x) = \rht(x) -\rhb(x)$, $\tilde{\beta} = \tfrac{\pi(y-\rhb(x))}{\tilde{h}(x)}$, we decompose $w(x,y)$ as in \eqref{eqn:w-decom}. Note that in the case of a flat top and bottom boundary, the coordinate change is trivial, and $w$ is identically equal to $v$. The case $y_0\geq\tfrac{1}{2}$ is treated in the analogous way by making a rotation about the top boundary.

We set $\tilde{S}(x) =  \tilde{S}^{\!B}(x)\cup \tS^{\!T}(x) = [\rhb(x),\rhb(x)+\tfrac{1}{4}]\cup [\rht(x)-\tfrac{1}{4},\rht(x)]$, and to establish the regularity of the nodal line, we first study it away from the boundary of $\tilde{\Omega}$. Define
\begin{equation}\label{e:e1}
e_1:=  \sup_{\substack{(x,y)\in \tilde{\Omega} \\ x \in I}}|\pa_x\tilde{E}(x,y)|  + \pi  \sup_{x \in I} |w_1(x)|\sup_{ x \in I}\left( \tih(x)^{-2}\left(|\tih'(x)| + |\rhb'(x)\rht(x)|\right)\right)
\end{equation}
and
\begin{equation}\label{e:lambda1}
\Lambda_1:=  \frac{1 }{2} \inf_{x\in  I}  \frac{|w_1'(x)|}{\tih(x)}   - e_1.
 \end{equation}
 We will show in the proof of Lemma \ref{lem:nodal-Lip} that $\Lambda_1$ provides a lower bound for $|\pa_xw(x_0,y_0)|$ for points $(\x,\y)\in w^{-1}(0)$ with $y_0\in\tilde{S}(x_0)^{c}$.

%%%%%%%%%%%%%%%%%%%%%%%%%%%%%%%%%%%%%%%%%%%%%%%%%%%%%
\begin{lem}\label{lem:nodal-Lip}
Suppose that $ \Lambda_1 >0$. Then, for  every $(\x,\y) \in w^{-1}(0)$ with  $\y \in \tS(\x)^c$ there exist a smooth real valued function $g$ and a neighborhood $U$ of $(\x,\y)$ such that
$$
w^{-1}(0) \cap U=\{(x,y)\in U: \;  x = g(y)\},
$$
with 
 \[|g'(y)| \leq\frac{1}{\Lambda_1}\Big(\pi \sup_{x \in I} \frac{|w_1(x)|}{\tih(x)} + \sup_{\substack{(x,y)\in \tilde{\Omega} \\ x \in I}}| \partial_y \tilde{E}(x,y)|\Big).\]

\end{lem}
%%%%%%%%%%%%%%%%%
\begin{proof1}{Lemma \ref{lem:nodal-Lip}} 
Note that for all $(x,y) \in \tilde{\Omega}$
\[
\pa_xw(x,y)  = w_1'(x)\sin(\tilde{\beta}(x,y))+w_1(x)\pa_x\tilde{\beta}(x,y)\cos(\tilde{\beta}(x,y)) + \pa_x\tilde{E}(x,y),
\] 
and 
\begin{equation}\label{e:beta der}
\pa_x\tilde{\beta}(x,y)=-\frac{\pi}{\tih(x)^2}\big( y\tih'(x)+\rhb'(x)\rht(x)\big).
\end{equation}
Therefore,
\[
|\pa_xw(x,y)|  \geq  |w_1'(x)|\sin(\tilde{\beta}(x,y)) -|w_1(x)\pa_x\tilde{\beta}(x,y)\cos(\tilde{\beta}(x,y))+\pa_x\tilde{E}(x,y)|.
\]
Let $(\x,\y) \in w^{-1}(0)$, and suppose that  $\y \in \tS(\x)^c$. Then, using $\sin(\tilde{\beta}(x_0,y_0))\geq \tfrac{1}{2\tilde{h}(x_0)}$
\begin{equation}\label{e:Lambda1*}
|\pa_xw(\x,\y)|  \geq   \frac{1}{2} \inf_{x\in  I}  \frac{|w_1'(x)|}{\tih(x)} - e_1=\Lambda_1,
\end{equation}
with $e_1$ as defined in \eqref{e:e1}.
Thus, $|\pa_xw(\x,\y)| >0$ since $\Lambda_1>0$ by assumption. This implies the existence of the graph function $g$ along a neighborhood of $w^{-1}(0)\cap \{(x,y)\in \Omega: \,y\in  \tS(x)^c\}$. Note that for every $y$  
\begin{equation}\label{e:g}
g'(y)=-\frac{\partial_y w(g(y),y)}{\partial_x w(g(y),y)}.
\end{equation}
We next find an upper bound for $|\partial_y w(g(y),y)|$.
Since for all $(x,y) \in \Omega$
\[
\partial_y w(x,y)=\pi \frac{w_1(x)}{\tih(x)} \cos(\tilde{\beta}(x,y)) + \partial_y \tilde{E}(x,y),
\]
then
\begin{equation} \label{e:y1-upper}
|\partial_y w(x,y)|\leq \pi \sup_{x \in I} \frac{|w_1(x)|}{\tih(x)} + \sup_{\substack{(x,y)\in \tilde{\Omega} \\ x \in I}}| \partial_y \tilde{E}(x,y)|.
\end{equation}
This together with \eqref{e:Lambda1*} yield the claimed bound on $|g'(y)|$ when $y \in \tS(x)^c$.
\end{proof1}

To study the regularity of the nodal line near $\pa\tilde{\Omega}$, we define
\begin{equation}\label{e:e2}
e_2:= \sup_{\substack{ x \in I\\ y \in  \tS^{\! \!B}(x)}} |\partial_y \partial_x \tilde{E}(x,y)|  + \pi \sup_{x \in I} |w_1(x)|\sup_{\substack{(x,y)\in \tilde{\Omega} \\ x \in I}} \frac{|\tilde{h}'(x)|}{\tih(x)^2}
\end{equation}
 and 
\begin{equation}\label{e:lambda2}
\Lambda_2:=  2 \inf_{x\in  I}  \frac{|w_1'(x)|}{\tih(x)}    -e_2.
 \end{equation}
We will show in Lemma \ref{lem:nodal-Lip2} that $\Lambda_2$ provides a lower bound for $|\pa_xw(x_0,y_0)|$ for all points $(\x,\y)\in w^{-1}(0)$.

\begin{lem} \label{lem:nodal-Lip2}
If $\Lambda_2>0$, there exist a  neighborhood $U$ of $(\x,\y)$, and smooth real valued function $g$,  such that
$
w^{-1}(0) \cap U=\{( x,  y)\in U: \;  x = g( y)\},
$
and
\[|g'( y)| \leq  \frac{| y|}{\Lambda_2} \, \sup_{\substack{ x \in I\\ y \in  \tS^{\!B}(x)}}\left( \tfrac{1}{2}  \tih(x)^{-2}(1+|\rhb''(x)|)  |\partial_y\tE(x,y)| + \tfrac{1}{2}  |\rhb''(x)|  |\partial_y^2E(x,y)| +  |\partial_y^3E(x,y)| \right)\]% \frac{1+ | y| \tfrac{1}{2}\sup_{x \in I}|\rhb(x)|}{\Lambda_2} \sup_{\substack{ x \in I\\ y \in  \tS^{\! \!B}(x)}} \Big[ | \partial_y^3\tE(x,y)||\rhb(x)|+ | \partial_y^3\tE(x,y)|+ 12 |\partial_y \tE(x,y)|  \Big] .\]
%In this system of coordinates, $|\tilde y| \leq y-\pb(x)$.
Furthermore, $w^{-1}(0)$ meets $\partial \tilde{\Omega}$ orthogonally.
\end{lem}
\begin{proof1}{Lemma \ref{lem:nodal-Lip2}}
%Next, assume  that $\y \in \S^{\! \!B}(\x)$.  To prove the existence of $g$ in a neighborhood of $(\x,\y)$ we will make a change coordinates $(x,y) \mapsto (\tilde x, \tilde y)$ and prove the existence of a function $\tilde g$ for which $v^{-1}(0)=\{(\tilde x, \tilde y): \; \tilde x=\tilde g(\tilde y)\}$  in a neighborhood of $(\x,\y)$.
%Let  $(\xx, \yy)\in \partial \Omega $ be so that $d( (\x,\y), \partial \Omega)=  d((\x,\y),(\xx,\yy))$. 
%Next,  let $F$ the linear isometry built by rotating about the point $(\x,\y)$ and then translating in such a way that if $(\x,\y)=F(\tilde \x, \tilde \y)$, then $(\xx,\yy)=F(\tilde \x,0)$. In general, write $(\tilde x, \tilde y)=F^{-1}(x,y)$ for the new system of coordinates. Note that by construction  $\tilde \y \leq \y-\pb(\x)$. 
%In addition, in these coordinates the bottom piece of $\partial \Omega$ consists of points of the form  $(\tilde x, \tilde \pb(\tilde x)) $ for some $ \tilde \pb$ with $\tilde \pb(\tilde \x)=0$ and $\tilde \pb'(\tilde \x)=0$.
%$\clubsuit$ CAUTION: In what follows I'm writing $(x,y)$ for these new coordinates $(\tilde x, \tilde y)$. The rest of the proof works if the solution $v$ is of the form $v(x,y) = \psi(x) \sin(\beta(x,y))+ E(x,y),$  {\bf where $(x, y)$ are these new coordinates. Here, $\psi, E,$ and $\beta$ (together with $\pb$ and $h$) are defined using these coordinates as well}. But it doesn't work if we only assume that $v$ has that form in the cartesian coordinate system. At the moment I don't know how to fix this in a clean way. $\clubsuit$\\
Since $\partial_{ x} \tE( x,  \rhb( x))=0$,  we have $\partial_{  x} \tE( \x,0)=0$, and so    
\[|\partial_{ x}  \tE( \x,  \y)| \leq  \y \sup_{y \in   \tS^{\!B}(\x)} |\partial_{ y} \partial_{x} \tE(\x,y)|.\]
  In addition, using  $\rhb( \x) = \rhb'( \x) =0$, we know that   $\sin(\tilde{\beta}( \x, \y))\geq \tfrac{2 }{\tih( \x)} \y$ and $\pa_x\tilde{\beta}( \x,  \y)=-\frac{\pi}{\tih( \x)^2} \y \tih'( \x)$.       Therefore, 
\begin{equation}\label{e:Lambda2*}
|\pa_xw(\x,\y)|  \geq   \y \Big( 2 \inf_{x\in  I}  \frac{|w_1'(x)|}{\tih(x)}    -e_2\Big)= \Lambda_2\,  \y.
\end{equation}
This proves the existence of $ g$. %The case $y \in \S^{\! \!T}(x)$ is treated in the same way.
 For $(\x,\y) \in w^{-1}(0)$, with $y_0\in \tS^{\!B}(\x)$ we have $w_1(\x)=-\frac{\tE(\x,\y)}{\sin (\tilde{\beta}(\x,\y))}$. Therefore, 
\[
\partial_y w(\x,\y) = -\frac{1}{\tih(\x)}\tE(\x,\y) \frac{\pi \cos (\tilde{\beta}(\x,\y))}{\sin (\tilde{\beta}(\x,\y))} + \partial_y \tE(\x,\y). 
\]

Note that  $\frac{\pi \cos(\pi s)}{\sin (\pi s)}=\tfrac{1}{s}(1 -r(s))$ with $0\leq r(s) \leq \tfrac{\pi^2}{2} |s|^2$ for $|s|<\tfrac{1}{2}$. 
Since $\tilde{\beta}(\x, \y)=\frac{\pi \y}{\tih(\x)}$ , it follows that 
\begin{align}\label{e:d_yv}
\partial_y w(\x,\y) 
&= -\tfrac{1}{\y}\tE(\x,\y) (1 - r(\tfrac{\y}{\tih(\x)}))+ \partial_y \tE(\x,\y) \notag\\
&= -\tfrac{1}{\y}\tE(\x,\y)+ \partial_y \tE(\x,\y)  +\tfrac{1}{\y}\tE(\x,\y)r(\tfrac{\y}{\tih(\x)}). 
\end{align}
Moreover, 
   \[
   \tfrac{1}{  \y}\tE( \x, \y)=\partial_y\tE( \x,0)+\tfrac{1}{2}\partial^2_y \tE( \x,0)\y+\tfrac{1}{6}\partial_y^3\tE( x,  y_1) \y^2,
   \] 
   for some $ y_1 \in \tS^{\!B}(\x)$, and
\[
\partial_y \tE( \x,  \y)=\partial_y \tE( \x,0)+\partial^2_y \tE( \x,0)\y+\tfrac{1}{2}\partial^3_y \tE( \x,  y_2)\y^2,
\]
for some $ y_2 \in \tS^{\!B}(\x)$.
In particular, \eqref{e:d_yv} yields
\begin{align}\label{e:d_yv 2}
\partial_y w(\x,\y) 
&=\tfrac{1}{2}\partial^2_y \tE( \x,0)\y+ \tfrac{1}{2}\partial^3_y \tE( \x,  y_2)\y^2 -\tfrac{1}{6}\partial_y^3\tE( x,  y_1) \y^2 +\tfrac{1}{\y}\tE(\x,\y)r(\tfrac{\y}{\tih(\x)}). 
\end{align}
 Next, note that since $\tE(\x, 0)=0$,  there exists $y_2 \in \tS^{\!B}(x) $ such that        $\tE(\x, \y)= \partial_y \tE(\x, y_3)\y$.
 Since $0\leq r(s) \leq \tfrac{\pi^2}{2} |s|^2$ for $|s|<\tfrac{1}{2}$, it follows that 
\begin{align*}
|\partial_y w(\x,\y)|
&\leq \tfrac{1}{2} |\partial^2_y \tE( \x,0)| \y+\tfrac{1}{2}|\partial_y^3\tE(\x,y_2)|\y^2+\tfrac{1}{6}|\partial_y^3\tE(\x,y_1)|\y^2+ \tfrac{\pi^2}{2}  \tih(\x)^{-2}|\partial_y \tE(\x,y_3)| \y^2.
\end{align*}
Therefore,
 \begin{equation}\label{e:Lambda2b beginner}
|\partial_y w(\x,\y)|
\leq    \tfrac{1}{2}\y|\partial^2_y \tE( \x,0)|+ \y^2 \sup_{\substack{ x \in I\\ y \in  \tS^{\!B}(x)}} \Big[| \partial_y^3\tE(x,y)|+ \tfrac{\pi^2}{2} \tih(\x)^{-2}|\partial_y \tE(x,y)|  \Big] .
\end{equation}
 In the same way, we have
\begin{align*}
    \pa_yw(x_0,0) = -\frac{\tfrac{\pi}{2\tih(\x)} \y^2\pa_y^2\tE(\x,y_1)}{\sin\left(\tfrac{\pi \y}{\tih(\x)}\right)} + \pa_y\tE(\x,0)r(\tfrac{\y}{\tih(\x)}), 
\end{align*}
and
\begin{align} \label{e:Lambda2b-0}
  \left|\pa_yw(x_0,0)\right| \leq \y \sup_{\substack{ x \in I\\ y \in  \tS^{\!B}(x)}} |\partial^2_y \tE( x,y)|+ \tfrac{\pi^2}{2} \tih(\x)^{-2}\y^2 |\partial_y \tE( \x,0)|.
\end{align}
To improve \eqref{e:Lambda2b beginner} and obtain a $\y^2$ in the upper bound, we need better control on  $\partial^2_y \tE( \x,0)$ in \eqref{e:d_yv 2}.  To do this we note
\begin{equation} \label{e:Lambda2bb}
\partial_y^2 \tE(\x,0)=\partial_y^2 w( \x,0)=-\partial_x^2 w(\x,0)-\mu w( \x,0)=-\partial_x^2 w(\x,0)=-\rhb''(\x)\partial_y w(\x, 0),
\end{equation}
where the last equality was obtained after differentiating $w(x, \rhb(x))\equiv 0$ twice and using  $\rhb'(\x)=\rhb(\x)=0$.
%Since $\partial_y v(\x,0)=\pi \frac{\psi(\x)}{h(x)} \cos(\beta(\x,0)) + \partial_y E(\x,0)= \pi \frac{\psi(\x)}{h(x)} +  \partial_y E(\x,0) ,$ we have
%$\partial^2_y v(x,y)=\pi \frac{\psi(x)}{h(x)} \cos(\beta(x,y)) + \partial_y E(x,y),$
%Next, note that since $0=v(x_0, y_0)=\psi(\x)\sin(\tfrac{\pi \y}{h(\x)})+E(\x,\y)$,
%\[
%|E(\x,\y)| \leq |\psi(\x)\sin(\tfrac{\pi \y}{h(\x)})| \leq   |\psi(\x)|  |\tfrac{\pi \y}{h(\x)}  \leq 2 \pi |\y| \sup_{x\in I}  |\psi(\x)|. 
%\]
From \eqref{e:Lambda2b beginner} it follows that 
for $y \in \tS^{\!B}(x)$ and $(x,y) \in w^{-1}(0)$,
\begin{equation}\label{e:Lambda2b}
%|\partial_y w(\x,\y)|
%\leq  \y^2(1+ \y \tfrac{1}{2}\sup_{x \in I}|\rhb(x)|) \sup_{\substack{ x \in I\\ y \in  \tS^{\! \!B}(x)}} \Big[ | \partial_y^3\tE(x,y)||\rhb(x)|+ | \partial_y^3\tE(x,y)|+ 12 |\partial_y \tE(x,y)|  \Big] .
|\partial_y w(\x,\y)|
\leq  \y^2\sup_{\substack{ x \in I\\ y \in  \tS^{\!B}(x)}}\left( \tfrac{\pi^2}{2} \tih(x)^{-2}(1+|\rhb''(x)|)  |\partial_y\tE(x,y)| + \tfrac{1}{2}  |\rhb''(x)|  |\partial_y^2E(x,y)| +  |\partial_y^3E(x,y)| \right)
\end{equation}
%{\color{red}I'm not sure why the estimate above follows from \eqref{e:Lambda2b beginner}. It seems like $\rhb''$ should appear, and we need something like \eqref{e:Lambda2b-0} as well as \eqref{e:Lambda2b beginner}.} 
The bound on $g'$ follows from combining \eqref{e:g} with \eqref{e:Lambda2*} and \eqref{e:Lambda2b}. In particular $g'(0)=0$, showing that $w^{-1}(0)$ meets $\pa\tilde{\Omega}$ orthogonally.
\end{proof1}
%%%%%%%%%%%%%%%%%%%%%%%%%%%%%%%%%
%\begin{rem} \label{rem:nodal-Lip2}
%Writing $|\pa_y^2w(x_0,y_0)| \leq |\pa_y^2w(x_0,0)| + y_0\sup_{y\in[0,y_0]}|\pa_y^3w(x_0,y)|$, and using the estimate in \eqref{e:Lambda2bb}, we see that $|\pa_y^2w(x_0,y_0)|$ is $O(y_0)$ as $y_0\to0$. We will use this in Lemma \ref{lem:second} below to obtain a uniform bound on $g''$ near the boundary. 
%\end{rem}
 
Applying Proposition \ref{prop:Main}, the estimate on $|g'(y)|$ given in Theorem \ref{thm:Main} follows immediately from Lemmas \ref{lem:nodal-Lip} and \ref{lem:nodal-Lip2}.  The following lemma gives the desired uniform bound on $g''(y)$ and completes the proof of Theorem \ref{thm:Main}:
\begin{lem} \label{lem:second}
There exist constants $c>0$, $C$ such that
\begin{align*}
\left|g''(y)\right| \leq C\left(\eta e^{-cN} + \frac{\delta}{N^2}\right).
\end{align*}
\end{lem}
\begin{proof1}{Lemma \ref{lem:second}}
Let $(x_0,y_0)$ be a point on the nodal line $\Gamma$. Then, differentiating $w(g(y),y) = 0$ twice gives the expression
\begin{align} \label{eqn:g-second}
 g''(y_0) = \frac{(\pa_yw(x_0,y_0))^2\pa_x^2w(x_0,y_0) + (\pa_xw(x_0,y_0))^2\pa_y^2w(x_0,y_0) -2 \pa_yw(x_0,y_0)\pa_xw(x_0,y_0)\pa_x\pa_yw(x_0,y_0)}{(\pa_xw(x_0,y_0))^3}.
 \end{align}
To bound the denominator, we use the lower bounds on $\pa_xw(x_0,y_0)$ from \eqref{e:Lambda1*} (in the centre) and \eqref{e:Lambda2*} (near the boundary). By Proposition \ref{prop:Main} this implies that
\begin{align*}
\left|\pa_xw(x_0,y_0)\right| \geq C^{-1}N^{-1}|y_0|.
\end{align*} 
We also have upper bounds on $\pa_yw(x_0,y_0)$ from \eqref{e:y1-upper} (in the centre) and \eqref{e:Lambda2b} (near the boundary), which again using Proposition \ref{prop:Main} gives
\begin{align*}
\left|\pa_yw(x_0,y_0)\right| \leq C|y_0|^2\left(\eta e^{-cN} + \frac{\delta}{N^3}\right).
\end{align*}
 Finally, from Proposition \ref{prop:Main} we have $\left|\pa_x\pa_yw(x_0,y_0)\right| \leq CN^{-1}$, $\left|\pa^2_xw(x_0,y_0)\right| \leq CN^{-2}$, and combining \eqref{e:Lambda2bb} with
\begin{align*}
\left|\pa^2_yw(x_0,y_0)\right| \leq \left|\pa^2_yw(x_0,0)\right| + |y_0|\sup_{y\in[0,y_0]}|\pa_y^3w(x_0,y_0)|,
\end{align*}
gives
\begin{align*}
\left|\pa_y^2w(x_0,y_0)\right| \leq C|y_0|\left(\eta e^{-cN} + \frac{\delta}{N^3}\right).
\end{align*}
Using these estimates gives the desired bound for the expression for $g''(y_0)$ in \eqref{eqn:g-second}.
\end{proof1}

\begin{rem} \label{rem:nodal-Lip}
In the case of flat upper and lower boundaries, we have the following estimates on the quantities appearing in the numerator of \eqref{eqn:g-second}: First, $\pa_yv(x_0,y_0)$ satisfies
\begin{align*}
    |\pa_yv(x_0,y_0)|& \leq \pi \sup_{x\in I}|v_1(x)| + \sup_{\substack{(x,y)\in {\Omega} \\ x \in I}}| \partial_y{E}(x,y)| \qquad \emph{ for } \tfrac{1}{4} \leq  y_0\leq \tfrac{3}{4}, \\
    |\pa_yv(x_0,y_0)|& \leq y_0^2  \sup_{\substack{(x,y)\in {\Omega} \\ x \in I}}\left[|\pa_y^3E(x,y) +  \tfrac{1}{2}  | \partial_y{E}(x,y)|\right] \qquad \emph{ for } y_0\leq \tfrac{1}{4}, y_0 \geq \tfrac{3}{4},
\end{align*}
where in the second inequality we have used \eqref{e:Lambda2b beginner} (and the fact that $\pa_y^2E(x_0,0) = 0$ in the flat case). We also immediately have the estimates on the second derivatives of $v$ of
\begin{align*}
|\pa_x\pa_yv(x_0,y_0)| & \leq \pi\sup_{x\in I}|v_1'(x)| + \sup_{\substack{(x,y)\in \Omega \\ x \in I}}|\pa_x\pa_yE(x,y)|, \\
|\pa_x^2v(x_0,y_0)| & \leq \sup_{x\in I}|v_1''(x)| + \sup_{\substack{(x,y)\in \Omega \\ x \in I}}|\pa_x^2E(x,y)|.
\end{align*}
Finally, 
\begin{align*}
    |\pa_y^2v(x_0,y_0)| & \leq \pi^2\sup_{x\in I}|v_1(x)| + \sup_{\substack{(x,y)\in \Omega \\ x \in I}}|\pa_y^2E(x,y)| \qquad \emph{ for } \tfrac{1}{4} \leq y_0\leq \tfrac{3}{4} \\
        |\pa_y^2v(x_0,y_0)| & \leq \pi^3 y_0\sup_{x\in I}|v_1(x)| +  y_0 \sup_{\substack{(x,y)\in \Omega \\ x \in I}}|\pa_y^3E(x,y)| \qquad \emph{ for } y_0\leq \tfrac{1}{4}, y_0 \geq \tfrac{3}{4},
\end{align*}
where in the second inequality, we have used $|\pa_y^2v(x_0,y_0)| \leq |\pa_y^2v(x_0,0)| + y_0\sup_{y\in[0,y_0]}|\pa_y^3v(x_0,y)|$, and that $\pa_y^2v(x_0,0)=0$ in the flat case. We will use these estimates in \emph{Section \ref{Fodes}} when we explicitly track the constants in the flat case.
\end{rem}

%%%%%%%%%%%%%%%%%%%%%%%%%%%%%%%%%%%%%
%%%%%%%%%%%%%%%%%%%%%%%%%%%%%%%%%%%%%
\section{Proof of Proposition \ref{prop:Main}}\label{sec: main prop}
%%%%%%%%%%%%%%%%%%%%%%%%%%%%%%%%%%%%%
%%%%%%%%%%%%%%%%%%%%%%%%%%%%%%%%%%%%%
\label{duhamel}
In this section we will prove Proposition \ref{prop:Main} by establishing the required properties of the decompositions of $v$ and $w$ defined in \eqref{eqn:v-decom} and \eqref{eqn:w-decom}. 
From the definition of the domain $\Omega$ from \eqref{eqn:Omega}, $\Omega$ contains the rectangle $[0,N]\times[\tfrac{\delta}{N^3},1-\tfrac{\delta}{N^3}]$, and is contained in the rectangle $[-2\eta,N+2\eta]\times[-\tfrac{\delta}{N^3},1+\tfrac{\delta}{N^3}]$. Therefore, by domain monotonicity for Dirichlet eigenvalues we have the following lemma.
\begin{lem}\label{lem:eigenvalue}
The second Dirichlet eigenvalue $\mu$ satisfies
\begin{align*}
\pi^2\left(1+2\tfrac{\delta}{N^3}\right)^{-2}  + 4\pi^2 (N+4\eta)^{-2}\leq \mu \leq \pi^2\left(1-2\tfrac{\delta}{N^3}\right)^{-2} +4\pi^2N^{-2}. 
\end{align*}
\end{lem}
The boundary of $\Omega$ is $C^2$ smooth, except at four points where the $C^2$-curves meet at a convex angle. This ensures that the gradient of $v(x,y)$ is bounded.
\begin{lem} \label{lem:gradient}
There exists a constant $C$ $($depending only on $\tilde{C}_1$, $\tilde{C}_2$ from \eqref{eqn:pb}$)$ such that $\left|\nabla v(x,y) \right| \leq C$ for $(x,y)\in\Omega$. In particular, from \eqref{eqn:pl}, for all $y$, we have 
\begin{align*}
\left|v(x,y)\right| \leq C(2\eta+x)\quad \emph{for } x\leq 1 , \qquad \left|v(x,y)\right| \leq C(2\eta+N-x) \quad\emph{for } x\geq N-1. 
\end{align*}
\end{lem}

We recall that the function $w$ satisfies $(\Delta + \mu)w = 0$ in the domain $\tilde{\Omega}$ as in \eqref{tOmegadef},
with $w(x,\rhb(x)) = w(x,\rht(x)) = 0$.  We recall that  $\rhb$, $\rht$, $\rhr$, $\rhl$, satisfy the bounds \eqref{e:rho B T} and \eqref{e:rho L R}. 
%\begin{align} \label{eqn:rhb}
%%\left|\rhb(x)\right| \leq \frac{2\delta}{N^3}(1+|x-\tilde{x}_0|), \quad \left|\rht(x) - 1\right|\leq \frac{2\delta}{N^3}(1+|x-\tilde{x}_0|), \quad %\left|\frac{d^j}{dx^j}\rhb(x)\right|, \left|\frac{d^j}{dx^j}\rht(x)\right| \leq 2\tilde{C}_j\frac{\delta}{N^3},
%\left|\frac{d^j}{dx^j}\rhb(x)\right|, \left|\frac{d^j}{dx^j}\rht(x)\right| \leq 2\tilde{C}_j\frac{\delta}{N^3},
%\end{align}
%$j\geq1$, and   $\rho_L$, $\rho_{R}$ satisfy
%\begin{align} \label{eqn:rhl}
%    -\eta- \frac{\delta}{N^3} \leq \rho_L({y})\leq \frac{\delta}{N^3}, \quad -\frac{\delta}{N^3} \leq \rho_R({y}) - N \leq \eta + \frac{\delta}{N^3}, \quad \left|\frac{d^j}{d{y}^j}\rho_L({y})\right|, \left|\frac{d^j}{d{y}^j}\rho_R({y})\right| \leq \eta + \frac{\delta}{N^3},
%\end{align}
Since $w$ is equal to $v$ in the rotated coordinates, by Lemma \ref{lem:gradient}  we have 
\begin{equation}
    \label{wbds}
\norm{w}_{L^{\infty}} = 1, |\nabla w(x,y)|\leq C.
\end{equation}
Defining a height function by $\tilde{h}(x) = \rhb(x) - \rht(x) \geq \tfrac{3}{5}$, we write $w$ as the Fourier series
\begin{align*}
w(x,y) = \sum_{k\geq1} w_k(x) \sin (k \tilde \beta(x,y)).
\end{align*}
 Here the $k$-th mode $w_k(x)$ is given by
\begin{align*}
w_k(x) = \frac{2}{\tilde{h}(x)} \int_{\rhb(x)}^{\rht(x)}w(x,y)\sin( k \tilde \beta(x,y)) \ud y.
\end{align*}
To prove Proposition \ref{prop:Main}, we will first bound each mode $w_k(x)$, then sum over $k$, and finally use elliptic estimates to extend these to derivative bounds. To estimate $w_k(x)$, we use the eigenfunction equation to find the equation that it satisfies, and then use the Duhamel principle to find an implicit expression. To bound this expression we need control on the boundary values $w_k(0)$, $w_k(N)$.
\begin{lem} \label{lem:boundary}
There exists a constant $C$ such that
\begin{align*}
|w_k(0)| + |w_k(N)| \leq C\left(\eta + \tfrac{\delta}{N^3}\right).
\end{align*}
\end{lem} 
\begin{proof1}{Lemma \ref{lem:boundary}}
By definition
\begin{align} \label{eqn:boundary-prop1}
w_k(0) =  \frac{2}{\tilde{h}(0)} \int_{\rhb(0)}^{\rht(0)}w(0,y)\sin( k \tilde \beta(0,y) ) \ud y.
\end{align}
As noted in Remark \ref{rem:rotation}, the angle of rotation in the definition of $w$ is bounded by $C\tfrac{\delta}{N^3}$ for some $C>0$.  Since $v$ vanishes on $\pa\Omega$, Lemma \ref{lem:gradient} implies that
\begin{align*}
|w(0,y)| \leq C\left(\eta + \tfrac{\delta}{N^3}\right).
\end{align*}
Inserting this into \eqref{eqn:boundary-prop1} gives the estimate for $|w_k(0)|$. The estimate for $w_k(N)$ follows in the same way. 
\end{proof1}

%%%%%%%%%%%%%%%%%%%%%%%%%%%%%%%%%%%%%
\subsection*{ Proof of Proposition \ref{prop:Main}: Flat case}
%%%%%%%%%%%%%%%%%%%%%%%%%%%%%%%%%%%%%
Let us first consider the case of a flat top and bottom, with $w=v$, and $\rhb(x) \equiv 0$, $\rht(x) \equiv 1$. In this case, we can remove the factor of $\tfrac{\delta}{N^3}$ in the estimate in Lemma \ref{lem:boundary} above. Using that $(\Delta+\mu)v(x,y) = 0$ for $0\leq x \leq N$, the function $v_k(x) = 2\int_{0}^{1}v(x,y)\sin(k\pi y)\ud y$ satisfies the ODE
\begin{align} \label{eqn:vk-ODE}
v_k''(x) + (\mu-\pi^2k^2)v_k(x) = 0.
\end{align}
Writing $\mu_k^2 = \pi^2k^2-\mu\geq \pi^2(k^2-1) - 4\pi^2N^{-2}\geq(k^2-2)\pi^2$ for $k\geq2$ (by Lemma \ref{lem:eigenvalue}, provided $N\geq 2$), we therefore have, for $k\geq 2$,
\begin{align}  \nonumber
v_k(x) & = \frac{1}{e^{\mu_kN}-e^{-\mu_kN}}\left(v_k(0)\left(e^{\mu_k(N-x)} - e^{\mu_k(x-N)}\right)  + v_k(N)\left(e^{\mu_kx} - e^{-\mu_kx}\right)\right) \\\label{eqn:vk}
& = \frac{1}{\sinh(\mu_kN)}\left(v_k(0)\sinh(\mu_k(N-x)) + v_k(N)\sinh(\mu_kx)\right). 
%v_k(x) = v_k(0)\left(1-e^{-2N\mu_k}\right)^{-1}\left(e^{-\mu_kx}-e^{\mu_k(x-2N)}\right).
\end{align}
Writing $v(x,y) = v_1(x)\sin(\pi y) + E(x,y)$ as in \eqref{eqn:v-decom}, this expression gives $\left|E(x,y)\right| \leq C\eta e^{-cN}$, and likewise for derivatives of $E(x,y)$. For $k=1$, \[ 4\pi^2(N+4\eta)^{-2} \leq \mu -\pi^2 \leq  4\pi^2N^{-2},\] and we set $\mu_1^2 = \mu-\pi^2$. The function $v_1(x)$ satisfies
\begin{align} \label{eqn:v1}
v_1(x) = v_1(0)\cos(\mu_1x) + \frac{v_1(N)-v_1(0)\cos(\mu_1N)}{\sin(\mu_1N)} \sin(\mu_1x).
%v_1(x) = A\sin(\mu_1(N-x)), \qquad v_1'(x) = \mu_1A\cos(\mu_1(N-x)), \qquad v_1''(x) = -\mu_1^2A\sin(\mu_1(N-x))
\end{align}
Setting $A_1 = \frac{v_1(N)-v_1(0)\cos(\mu_1N)}{\sin(\mu_1N)}$ gives $|v_1(x) - A_1\sin(\mu_1x)| \leq C\eta$, and since $\norm{v}_{L^{\infty}} = 1$, this implies that $\left||A_1|-1\right|\leq C\eta$. The estimates from Proposition \ref{prop:Main} then follow readily from the expressions in \eqref{eqn:vk} and \eqref{eqn:v1}. %Therefore, $C^{-1}N^{-1} \leq |v_1'(x)| \leq CN^{-1}$, and we have proved Proposition \ref{prop:Main} in this flat case.

%%%%%%%%%%%%%%%%%%%%%%%%%%%%%%%%%%%%%
\subsection*{Proof of Proposition \ref{prop:Main}: General case}
%%%%%%%%%%%%%%%%%%%%%%%%%%%%%%%%%%%%%

In the general case, $w_k(x)$ satisfies an approximate version of the ODE in \eqref{eqn:vk-ODE}, with an error depending on $\rhb(x)$, $\rht(x)$ and their first two derivatives: Fix $x^*\in[0,N]$, and set 
\begin{equation}
    \label{ekdef}
e_k(x,y) = \frac{2}{\tilde{h}(x)}\sin( k \tilde \beta(x,y)).
\end{equation}

\begin{lem} \label{lem:vk-eqn}
The function $w_k(x)$ satisfies the equation
\begin{align*}
w_k''(x) + \left(\mu - \frac{\pi^2k^2}{\tilde{h}(x^*)^2} \right)w_k(x) = F_k(x),
\end{align*}
where $F_k(x)$ has the bound
\begin{align*}
|F_k(x)| \leq Ck\left(\left| \tfrac{1}{\tilde{h}(x)^2}-\tfrac{1}{\tilde{h}(x^*)^2}\right| + \left|\rht'(x)\right|+\left|\rhb'(x)\right| + \left|\rht''(x)\right|+\left|\rhb''(x)\right|\right),
\end{align*}
for an absolute constant $C$.% In particular, at $x=x^*$,  $F_k(x^*)$ can be bounded by terms only depending linearly on $k$. 
\end{lem}
\begin{proof1}{Lemma \ref{lem:vk-eqn}}
The function $F_k(x)$ is equal to
\begin{align} \label{eqn:Fk-exp}
F_k(x) = \pi^2k^2\left(\frac{1}{\tilde{h}(x)^2}-\frac{1}{\tilde{h}(x^*)^2}\right)w_k(x) + 2\int_{\rhb(x)}^{\rht(x)} \pa_xw(x,y)\pa_xe_k(x,y) \ud y + \int_{\rhb(x)}^{\rht(x)}w(x,y) \pa_x^2e_k(x,y)\ud y.
\end{align}
%Since we will need to sum over $k$, we need to consider the dependence on $k$, as well as the size and regularity of each term. In particular, the highest power of $k$ appearing in $F_k(x)$ is $k^2$, and the highest derivative of $\phi_T(x)$ appearing in $F_k(x)$ is $\phi_T''(x)$. Moreover, no term in $F_k(x)$ contains both a factor of $k^2$ and $\phi_T''(x)$. 
Applying the bounds we have derived for $w$ in \eqref{wbds} and the definition of $e_k$ in \eqref{ekdef}, the terms that do not immediately obey the estimate of the lemma are the first term and the term in the final integral given by
\begin{align} \label{eqn:Fk-exp1}
-\frac{2k^2}{\tilde{h}(x)}  \int_{\rhb(x)}^{\rht(x)}w(x,y) \left(\pa_x\tilde{\beta}(x,y)\right)^2 \sin\left(k \tilde{\beta}(x,y) \right) \ud y.
\end{align}
This is because this is the only term in the final integral in \eqref{eqn:Fk-exp} for which a factor $k^2$ appears in the expression for $\pa_x^2e(x,y)$. All of the other terms in the last two integrals in \eqref{eqn:Fk-exp} contain at most one derivative of $w$, two derivatives of $\rht$ and $\rhb$, and one factor of $k$.  After an integration by parts $w_k(x)$ is equal to
\begin{align*}
\frac{2}{k\pi}\int_{\rhb(x)}^{\rht(x)} \pa_y w(x,y) \cos\left(k\tilde{\beta}(x,y)\right) \ud y,
\end{align*}
 and \eqref{eqn:Fk-exp1} can be written as
\begin{align*}
-\frac{2k}{\pi} \int_{\rhb(x)}^{\rht(x)}\pa_y \left( w(x,y) \left(\pa_x\tilde{\beta}(x,y)\right)^2 \right) \cos\left(k\tilde{\beta}(x,y)\right) \ud y .
\end{align*}
Since $|\pa_yw(x,y)|$ is bounded by a constant, both of these terms are of the desired form. 
\end{proof1}
The function $w_k(x)$ also satisfies the boundary conditions
\begin{align*}
w_k(0) = \alpha^{(1)}_k, \qquad w_k(N) = \alpha^{(2)}_k,
\end{align*}
where $\alpha^{(i)}_k$ are values coming from the side variation of the domain, wih bounds in Lemma \ref{lem:boundary}.
\\
\\
For $k=1$, set $\mu_1^2 = \mu-\tfrac{\pi^2}{\tilde{h}(x^*)^2}\geq0$ and for $k\geq 2$, set $\mu_k^2 = \mu_k(x^*)^2 = \tfrac{\pi^2 k^2}{\tilde{h}(x^*)^2}-\mu\geq0$.
\begin{lem} \label{lem:v1-vk}
Define the functions $W_1(x)$ and $W_k(x)$ (for $k\geq2$) by 
\begin{align*}
W_1(x) = \frac{1}{\mu_1}\sin(\mu_1 x), \qquad W_k(x) = \frac{1}{2\mu_k}\left(e^{\mu_kx} - e^{-\mu_kx}\right) = \frac{1}{\mu_k}\sinh(\mu_kx).
\end{align*}
Then,
\begin{align*}
w_1(x) = \int_{x}^{N}W_1(t-x)F_1(t)\ud t + A_1\sin(\mu_1 (N-x)) + \alpha_1^{(2)}\cos(\mu_1(N-x)),
\end{align*}
 with $\alpha^{(1)}_1 = \int_0^NW_1(t)F_1(t) \ud t + A_1\sin(\mu_1N) + \alpha_1^{(2)}\cos(\mu_1N)$. Also,
\begin{align*}
w_k(x) = \int_0^x W_k(x-t)F_k(t) \ud t + A_ke^{\mu_kx} + B_ke^{-\mu_kx},
\end{align*}
for constants $A_k$, $B_k$, with
\begin{align*}
\alpha^{(1)}_k & = A_k+B_k \\
\alpha^{(2)}_k & = A_ke^{\mu_kN}+B_ke^{-\mu_kN} + \int_0^NW_k(N-t)F_k(t)\ud t.
\end{align*}
\end{lem}
\begin{proof1}{Lemma \ref{lem:v1-vk}}
The functions $W_1(x)$ and $W_k(x)$ satisfy
\begin{align*}
W_1''(x) + \mu_1^2W_1(x) & = 0, \qquad W_1(0) = 0, \quad W_1'(0) = 1 \\
W_k''(x) - \mu_k^2W_k(x) & = 0, \qquad W_k(0) = 0, \quad W_k'(0) = 1. 
\end{align*}
The lemma then follows from Lemma \ref{lem:vk-eqn} and the boundary conditions of $w_k(x)$ at $x=0,N$.
\end{proof1}
Combining Lemmas \ref{lem:vk-eqn} and \ref{lem:v1-vk}, we can bound $w_k(x)$.
\begin{prop} \label{prop:wk-bound}
There exist constants $c$, $C$, such that for $x\in[0,N]$, $k\geq2$,
\begin{align*}
|w_k(x)|\leq C\left(\eta e^{-c\mu_k d(x)} + k^{-1}\frac{\delta}{N^3}\right). 
\end{align*}
Here $d(x) = \min\{x,N-x\}$ is the distance of $x$ from the endpoints of $[0,N]$.
%{\color{red} This is assuming that $|w_k(0)|, |w_k(N)| \leq C\left(\eta + \tfrac{\delta}{N^3}\right)$.}
\end{prop}
\begin{proof1}{Proposition \ref{prop:wk-bound}}
We fix $x^*\in[0,N]$, and use Lemmas \ref{lem:vk-eqn} and \ref{lem:v1-vk} to bound $w_k$ at $x=x^*$ (with a bound independent of $x^*$). 
The constants $A_k$, $B_k$ from Lemma \ref{lem:v1-vk} can be written for $k\geq2$ as
\begin{align*}
A_k = -\frac{1}{2\mu_k}\int_0^N\frac{W_k(N-t)}{W_k(N)}F_k(t) \ud t + \frac{-e^{-\mu_kN}\alpha^{(1)}_k+\alpha^{(2)}_k}{e^{\mu_kN}-e^{-\mu_kN}},
\end{align*}
with $B_k = \alpha^{(1)}_k-A_k$. From Lemma \ref{lem:vk-eqn}, we have the bound 
\begin{align} \label{eqn:wk-bound2}
|F_k(t)| \leq Ck \frac{\delta}{N^3}(1+|t-x^*|), 
\end{align}
and from Lemma \ref{lem:boundary},  $|w_k(0)|, |w_k(N)| \leq C\left(\eta + \tfrac{\delta}{N^3}\right)$. Therefore, since $\mu_k\geq\pi \sqrt{k^2-2}$, the only terms in the expression for $w_k(x^*)$ from Lemma \ref{lem:v1-vk} that do not immediately satisfy the required estimates are
\begin{align*} 
& \int_{0}^{x^*} W_k(x^*-t)F_k(t)\ud t - \frac{1}{2\mu_k}e^{\mu_kx^*}\int_{0}^{N} \frac{e^{\mu_k(N-t)}}{e^{\mu_kN}}F_k(t)\ud t .
\end{align*}
However, these integrals can be combined to be written as
\begin{align} \label{eqn:wk-bound1}
 -\frac{1}{2\mu_k}\int_{0}^{x^*} e^{-\mu_k(x^*-t)}F_k(t)\ud t  - \frac{1}{2\mu_{k}}\int_{x^*}^{N}e^{\mu_k(x^*-t)}F_k(t)\ud t.  
\end{align}
Using the bound on $F_k(t)$ from \eqref{eqn:wk-bound2} and integrating gives the desired bound.
\end{proof1}
We write
\begin{align*}
w(x,y) = w_1(x)\sin\left(\tilde{\beta}(x,y)\right) + \sum_{k\geq2}w(x,y)\sin\left(k\tilde{\beta}(x,y) \right)  = V_1(x,y) + \tilde{E}(x,y).
\end{align*}
Summing the estimate from Proposition \ref{prop:wk-bound} over $k$ we can control the $L^2$-norm of $\tilde{E}$. For the rest of the section, fix $x^*\in[1,N-1]$ and denote the cross-section at $x^*$ by $U(x^*) = \tilde{\Omega}\cap\{(x,y):x=x^*\}$.
\begin{cor} \label{cor:wk-bound}
 There exist constants $c$, $C$ such that
\begin{align*}
\|\tilde{E}\|_{L^2(U(x^*))} \leq C\left(\eta e^{-cd(x^*)} + \frac{\delta}{N^3}\right).
\end{align*}
\end{cor}
We now convert this $L^2$-estimate into bounds on derivatives of $\tilde{E}$.
\begin{prop} \label{prop:tildeE-bound}
For each $j\geq0$, and with $c>0$ as in Corollary \ref{cor:wk-bound}, there exists a constant $C_j$ such that 
\begin{align*}
\|{\tilde{E}}\|_{H^{j}(U(x^*))} \leq C_j\left(\eta e^{-c d(x^*)}  + \frac{\delta}{N^3}\right). 
\end{align*}
\end{prop}
\begin{proof1}{Proposition \ref{prop:tildeE-bound}}
To obtain this estimate on $\tilde{E}$ we find the elliptic equation that it satisfies. For $V_1(x,y):=w_1(x)\sin\left(\frac{\pi (y-\rhb(x))}{\tilde{h}(x)}\right)$, we have
\begin{align*}
\Delta V_1(x,y) & = \frac{2}{\tilde{h}(x)} \left(\int_{\rhb(x)}^{\rht(x)}\pa_x^2w(x,y') \sin \left(\tilde{\beta}(x,y') \right) \ud y' \right) \sin \left(\tilde{\beta}(x,y) \right) - \frac{\pi^2}{\tilde{h}(x)^2}V_1(x,y) + G_1(x,y) \\
& = -\mu V_1(x,y) + G_1(x,y).
\end{align*}
The function $G_1(x,y)$ consists of terms where at least one derivative in $x$ has been applied to a factor of $\rht(x)$ or $\rhb(x)$, and so for each $j\geq0$, there exists a constant $C_j$ such that 
\begin{align} \label{eqn:G1-bound}
\norm{G_1}_{H^{j}(U(x^*))} \leq C_j\frac{\delta}{N^3}. 
\end{align}
Using the eigenfunction equation, $\tilde{E}(x,y)$ satisfies
\begin{eqnarray*}
    \left\{ \begin{array}{rlc}
    \Delta \tilde{E}(x,y) & = -\mu \tilde{E}(x,y) - G_1(x,y) & \text{in } \tilde{\Om} \\
   \tilde{E}(x,y) & = 0& \text{on } \pa \tilde{\Om}.
    \end{array} \right.
\end{eqnarray*}
Applying elliptic estimates to this equation,  \eqref{eqn:G1-bound} and the estimate on $\tilde{E}$ from Corollary \ref{cor:wk-bound} establishes the proposition.
\end{proof1}

Using  Proposition \ref{prop:tildeE-bound}, we can obtain more refined information about the first Fourier mode $w_1(x)$ and complete the proof of Proposition \ref{prop:Main}.
\begin{prop} \label{prop:w1-control}
There exists a constant $C$ such that in the interval $[\tfrac{N}{4},\tfrac{3N}{4}]$ the function $w_1(x)$ has a unique zero at $x=x_0$ with $\left|x_0-\tfrac{N}{2}\right| \leq CN(\eta + \delta)$. Moreover, $|w_1'(x)|\geq C^{-1}N^{-1}$ for this range of $x$, and for $x\in[0,N]$, $1\leq j \leq 3$, we have,
\begin{align*}
\left|w_1(x) - A_1\sin(\mu_1(N-x))\right| \leq C\left( \eta  +\delta \right),  \qquad |w^{(j)}_1(x)| \leq CN^{-j}
\end{align*}
Here the constant $A_1$ is as in \emph{Lemma \ref{lem:v1-vk}} and satisfies $\left||A_1|-1\right|\leq C(\eta  + \delta)$.
\end{prop}
\begin{proof1}{Proposition \ref{prop:w1-control}}
By Lemma \ref{lem:eigenvalue}, $\left|\mu_1 - \tfrac{2\pi}{N}\right| \leq CN^{-2}\left(\delta+\eta\right)$ for a constant $C$. Therefore, using the expression for $w_1(x)$ from Lemma \ref{lem:v1-vk} and the bound $\left|F_1(t)\right| \leq C\frac{\delta}{N^3}(1+|t-x^*|)$ from Lemma \ref{lem:vk-eqn}, we have
\begin{align} \label{eqn:w1-control1}
\left|w_1(x) - A_1\sin(\mu_1(N-x))\right| \leq C(\eta  + \delta ).
\end{align}
Here $C$ is a constant (changing from line-to-line). Moreover, since $\norm{w}_{L^{\infty}} = 1$, and \[w(x,y) = w_1(x) \sin\left(\tilde{\beta}(x,y) \right) + \tilde{E}(x,y),\] combining \eqref{eqn:w1-control1} with Proposition \ref{prop:tildeE-bound}, we have
\begin{align*}
\left||A_1|-1\right|\leq C( \eta  +\delta).
\end{align*}
To complete the proof of the lemma, we need to bound $w_1'(x)$. Differentiating the expression from Lemma \ref{lem:v1-vk} gives
\begin{align*}
w_1'(x) = - \int_{x}^{N}W_1'(t-x)F_1(t)\ud t - \mu_1A_1\cos(\mu_1(N-x)) +\mu_1\alpha^{(2)}_1\sin(\mu_1(N-x)).
\end{align*}
In particular $\left|w_1'(x) +  \mu_1A_1\cos(\mu_1(N-x))\right| \leq CN^{-1}(\eta+ \delta)$, and combining this with the estimate for $A_1$ gives the required bound for $|w_1'(x)|$. The expression for $w_1''(x)$ from Lemma \ref{lem:vk-eqn} gives $|w_1''(x)|\leq CN^{-2}$, and differentiating we have $|w_1'''(x)|\leq C N^{-3}$. Since $|w_1'(x)|$ is non-zero on $[\tfrac{N}{4},\tfrac{3N}{4}]$, $w_1(x)$ has at most one zero in this interval. The function $\sin(\mu_1(N-x))$ has its unique zero in this interval at $\tilde{x}_0$, with $|\tilde{x}_0-\tfrac{N}{2}|\leq C(\delta+\eta)$, and its derivative is bounded below by $C^{-1}N^{-1}$. Therefore, $w_1(x)$ also has a unique zero at $x=x_0$, with $|x_0-\tfrac{N}{2}|\leq CN(\eta  + \delta)$.
\end{proof1}
\begin{rem} \label{rem:w1-control}
In the case that no rotation has been applied (so that $v_1 = w_1$), the function $F_1(t)$ from Lemma \ref{lem:vk-eqn} satisfies the stronger bound $|F_1(t)| \leq C\frac{\delta}{N^3}$. Inserting this stronger estimate into the argument above, the point $x_0$ satisfies  $|x_0-\tfrac{N}{2}|\leq CN(\eta  + \delta N^{-1})$.
\end{rem}

%%%%%%%%%%%%%%%%%%%%%%%%%%%%%%%%%%%%%
%%%%%%%%%%%%%%%%%%%%%%%%%%%%%%%%%%%%%
\section{An explicit Hadamard variation formula and constant tracking}
\label{Fodes}
%%%%%%%%%%%%%%%%%%%%%%%%%%%%%%%%%%%%%
%%%%%%%%%%%%%%%%%%%%%%%%%%%%%%%%%%%%%

 To prove Proposition \ref{prop:Main}, we used the estimate on the boundary values of the Fourier modes, $w_k(0)$ and $w_k(N)$, in Lemma \ref{lem:boundary}, which follows only from a gradient estimate on the eigenfunction. In order to track the constants appearing in the error estimates in Proposition \ref{prop:Main} in the flat case, we require a more explicit bound on these boundary values. We do this as follows, using a variant of a calculation given in \cite{GJ-Pacific}. 

\begin{prop} \label{prop:Hadamard-flat}
There exists a constant $C$ such that 
\begin{align} \label{eqn:vk-main}
\left|w_k(0) - \frac{4\pi}{N} \int_{\rhb(0)}^{\rht(0)}\rho_L(y)\sin(\tilde{\beta}(0,y))\sin(k\tilde{\beta}(0,y) ) \ud y \right|\leq  C\eta\left(\eta^{3/4} + \delta  +k^2(\eta+\delta N^{-3})^3\right).
\end{align}
\end{prop}
\begin{proof1}{Proposition \ref{prop:Hadamard-flat}}
%Let us assume that for the class of perturbations $\phi(y)$ under consideration, there exists an absolute constant $C_1$ such that
%\begin{align} \label{eqn:v-bound}
%\left|v(x,y)\right| \leq C_1(\eta+|x|) \qquad \left|\nabla v(x,y)\right| \leq C_1 \text{ for } x\leq 1.
%\end{align}
%Under this assumption, we can show that there exists a constant $C$ such that 
%\begin{align} \label{eqn:vk-main}
%\left|v_k(0) + \frac{4\pi}{N}\int_{0}^{1}\phi(y)\sin(\pi y)\sin(k\pi y) \ud y\right|\leq  C\left(\eta^{7/4}+k^2\eta^3\right).
%\end{align}
%Using the bounds from Lemma \ref{lem:gradient}, there exists a constant $C_1$ such that 
%\begin{align*}
%|v_k(0)| = 2\left|\int_{0}^{1}v(0,y)\sin(k\pi y)\ud y\right| \leq C_1\eta,
%\end{align*}
%and likewise for $v_k(N)$. 
%In the final step, we have used $|\mu_1 -\tfrac{2\pi}{N}| \leq C \tfrac{\eta}{N^2}$. %Combining \eqref{eqn:E-first} and \eqref{eqn:E-second}, we obtain
%\begin{align} \label{eqn:E-third}
%|E(x,y)|\leq C_4\eta^{1/2}.
%\end{align}
%This also implies the lower bound on $A$ of
%\begin{align} \label{eqn:A-lower}
%A \geq 1-C_5\eta e^{-C_5N}.
%\end{align}
We extract the main term in $w_k(0)$ as follows. First, we integrate by parts to write $w_k(0)$ as
\begin{align}\nonumber
2\int_{\rhb(0)}^{\rht(0)}& w(0,y) \sin(k\tilde{\beta}(0,y)) \ud y = 2\int_{\pa\tilde{\Omega}_0}w(x,y)\frac{\pa}{\pa \nu}\left(x\sin(k\tilde{\beta}(0,y))\right) \ud \sigma \\ \label{eqn:vk0}
& = -2\int_{\pa\tilde{\Omega}_0}\frac{\pa w}{\pa\nu}(x,y) x\sin(k\tilde{\beta}(0,y))\ud \sigma + 2\left(\mu-\frac{k^2\pi^2}{\tilde{h}(0)^2}\right)\int_{\tilde{\Omega}_0}w(x,y) x\sin(k\tilde{\beta}(0,y)) \ud x \ud y %\\ %\label{eqn:vk0}
%& = -2\int_{0}^{1} \left(\pa_x-\pl'(y)\pa_y\right)v(\pl(y),y) \pl(y)\sin(k\pi y) \ud y  + 2\left(\mu-\frac{k^2\pi^2}{\tilde{h}(0)^2}\right)\int_{\tilde{\Omega}_0}v(x,y) x\sin(k\beta(0,y)) \ud x \ud y . 
\end{align}
The domain $\tilde{\Omega}_0$ is the domain $\{(x,y)\in\tilde{\Omega}:\rho_L(y)\leq x \leq0\}$. The second integral in \eqref{eqn:vk0} is bounded in absolute value by $C k^2(\eta+\delta N^{-3})^3$. The first integral in \eqref{eqn:vk0} consists of three terms. Two of these integrals are over portions of the top and bottom boundaries of $\tilde{\Omega}$ of length bounded by $C(\eta+\delta N^{-3})$, and so since the gradient of $w$ is bounded, these integrals are bounded in absolute value by $C (\eta+\delta N^{-3})^2$. The remaining contribution to \eqref{eqn:vk0} is given by
\begin{align} \label{eqn:vk0a}
-2\int_{\rhb(0)}^{\rht(0)} \left(\pa_x-\rho_L'(y)\pa_y\right)w(\rho_L(y),y) \rho_L(y)\sin(k \tilde{\beta}(0,y) ) \ud y.
\end{align}
To pick out the main term in \eqref{eqn:vk0a} we write
\begin{align*}
    w(x,y) = \sin(\mu_1(N-x))\sin( \tilde{\beta}(0,y) ) + B(x,y),
\end{align*}
with $\mu_1^2 = \mu_1^2(0) = \mu-\tfrac{\pi^2}{\tih(0)^2}$ and for an error term $B(x,y)$ to be estimated below. Using $|\mu_1-\tfrac{2\pi}{N}| \leq CN^{-2}(\eta +\delta)$, \eqref{eqn:vk0a} becomes
\begin{align} \label{eqn:vk0b}
\frac{4\pi}{N} \int_{\rhb(0)}^{\rht(0)} \ \rho_L(y)\sin(\tilde{\beta}(0,y)) \sin(k\tilde{\beta}(0,y)) \ud y  -2\int_{\rhb(0)}^{\rht(0)} \left(\pa_x-\rho_L'(y)\pa_y\right)B(\rho_L(y),y) \rho_L(y)\sin(k\tilde{\beta}(0,y)) \ud y
\end{align}
up to an error $C(\eta+\delta)(\eta + \delta N^{-3})$. We are left to bound the second integral in \eqref{eqn:vk0b}, and to do this we will use the results of Section \ref{duhamel} to estimate $B(x,y)$. Summing the estimate from Proposition \ref{prop:wk-bound} over $k\geq2$, we obtain a bound on $\tilde{E}(x,y)$ for $x\geq0$ of
\begin{align*} 
\norm{\tilde{E}(x,y)}_{L^2(U(x))} \leq C\left(\frac{\eta}{\max\{1,x\}}e^{-c x} + \frac{\delta}{N^3}\right),
\end{align*}
for constants $c$, $C$, where we recall that $U(x)$ is the cross-section of $\tilde{\Omega}$ at $x$. Combining this with Proposition \ref{prop:w1-control} shows that for $0\leq x \leq 1$, we have
\begin{align}\label{eqn:B-first}
    \norm{B(x,y)}_{L^2(U(x))} \leq C\left(\frac{\eta}{\max\{1,x\}}e^{-c x} + \frac{\delta}{N^3}\right) + C(\eta+ \delta)
\end{align}
%Since we have normalised $v$ to have $L^{\infty}$-norm equal to $1$, \eqref{eqn:E-first} implies the existence of a constant $c_1>0$ so that
%\begin{align} \label{eqn:A-upper}
%A \leq 1+C_1\eta e^{-c_1N}.
%\end{align}
Using Lemma \ref{lem:gradient}, we can also bound $B(x,y)$ in a different way for $0\leq x\leq1$ via
\begin{align} \label{eqn:B-second}
|B(x,y)| \leq |w(x,y)| + \left|\sin(\mu_1(N-x))\sin( \tilde{\beta}(0,y) )\right| \leq C(\eta +\delta N^{-3}+x).% + \left(1+\eta\right)\sin(\mu_1x)) \leq C\left(\eta + \tfrac{\delta}{N^3} + x\right).
\end{align}
In particular, using \eqref{eqn:B-first} and \eqref{eqn:B-second}, we have  $\norm{B}_{L^2(\tilde{\Omega}_1)}\leq C\left(\eta^{3/4}+\delta\right)$. Moreover, $B$ satisfies the equation
\begin{align*}
    \Delta B = \Delta w + \left(\mu_1^2 + \tfrac{\pi^2}{\tih(0)^2} \right)(w - B) = - \left(\mu_1^2 + \tfrac{\pi^2}{\tih(0)^2} \right)B. 
\end{align*}
We can use this to bound the second integral in \eqref{eqn:vk0b}. Let $\chi(x)$ be a smooth cut-off function, equal to $1$ for $x\leq \tfrac{1}{4}$ and $0$ for $x\geq \tfrac{3}{4}$. There exists an extension $H(x,y)$ of $\sin(\mu_1(N-x))\sin( \tilde{\beta}(0,y) )\big|_{x=\rho_L(y)}$ to $\tilde{\Omega}_1$, with $H(1,y)\equiv0$ such that
\begin{align*}
\norm{H}_{H^{1}(\Omega_1)}  \leq C\eta.
\end{align*}
The function $\chi(x)B(x,y) - H(x,y)$ therefore vanishes on $\pa\tilde{\Omega}_1$, and satisfies
\begin{align*}
\Delta \left(\chi(x)B(x,y)-H(x,y)\right) = (\Delta \chi)B + 2\nabla\chi.\nabla B + \chi \Delta B - \Delta H =: F.
\end{align*}
Elliptic estimates therefore imply that
\begin{align*}
\norm{\chi B-H}_{H^{1}(\tilde{\Omega}_1)} \leq \norm{F}_{H^{-1}(\tilde{\Omega}_1)} \leq C\left(\norm{B}_{L^2(\tilde{\Omega}_1)} + \norm{H}_{H^{1}(\tilde{\Omega}_1)}\right).
\end{align*}
Using $\norm{B}_{L^2(\tilde{\Omega}_1)}\leq C(\eta^{3/4}+\delta )$, we have the same bound on $\norm{\chi B-H}_{H^{1}(\Omega_1)}$. Elliptic estimates thus give $\norm{\frac{\pa B}{\pa \nu}}_{L^2(D)} \leq C(\eta^{3/4}+\delta)$, where $D = \{(\rho_L(y),y):\rhb(0)\leq y \leq \rht(0)\}$. Applying this estimate  in the second integral in \eqref{eqn:vk0b}, we see that the integral can be bounded by $C(\eta+\delta N^{-3})(\eta^{3/4}+\delta )$, and this therefore concludes the proof of the proposition.
\end{proof1}

Now consider the case where the domain $\Omega$ has flat top and bottom boundaries (so that $\pt(x)=1$, $\pb(x)=0$). Proposition \ref{prop:Hadamard-flat} then allows us to track the constants appearing in the error estimates in Proposition \ref{prop:Main}: Given $N$ (not necessarily large) and a small constant $c>0$, we can choose $\eta = \eta(N)$ sufficiently small so that
\begin{align*}
|v_k(0)|, |v_k(N)| \leq \frac{8 \eta}{N} + c\left(\frac{\eta}{N} + \frac{k^2\eta}{N^2}\right).
\end{align*}
By choosing $c$ small compared to $8$, we can use this in \eqref{eqn:vk} and \eqref{eqn:v1} to get explicit estimates on $E(x,y)$ and its derivatives, with the estimates not depending on any unknown constants. Using this in the quantities appearing in Section \ref{sec:nodal}, for any $N \geq 8$ fixed and $\eta = \eta(N)$ sufficiently small, this provides the following bounds and proves Corollary \ref{cor:Main}: We have
\begin{align*}
    \tau \leq 10^{-4}\eta,
\end{align*}
where we recall from Lemma \ref{lem:interval2} that the width of the nodal line is bounded by $2\tau$. Moreover, we have
\begin{align*}
    \Lambda_1 \geq 0.3, \qquad \Lambda_2 \geq 1.2,
\end{align*}
and from Lemmas \ref{lem:nodal-Lip} and \ref{lem:nodal-Lip2} this gives the upper bound on $|g'(y)|$ of
\begin{align*}
    |g'(y)|\leq 10^{-2}\eta.
\end{align*}
Finally, inserting the bounds from Remark \ref{rem:nodal-Lip} on the terms appearing in $g''(y)$ gives
\begin{align*}
    |g''(y)|\leq 10^{-2}\eta.
\end{align*}

\bibliographystyle{alpha}
\bibliography{Refs}

\end{document}